\documentclass[a4paper,11pt]{amsart}

\usepackage{amsmath,amssymb,amsthm,mathtools,color,fullpage,graphicx,hyperref}
\usepackage[english]{babel}
\usepackage{stmaryrd}
\usepackage{dsfont}
\usepackage{bbold}
\usepackage[all]{xy}
\usepackage{csvsimple}
\usepackage{nicefrac}

\newcommand{\e}{\varepsilon}
\newcommand{\field}{\mathbb{k}}
\newcommand{\R}{\mathbb{R}}

\newcommand{\N}{\mathbb{N}}

\newcommand{\dist}{\mathrm {d}}
\newcommand{\distb}{\dist_{\mathrm{b}}}

\newcommand{\distone}{\dist_1}
\newcommand{\disti}{\dist_{\mathrm{i}}}

\newcommand{\hdistb}{\hat\dist_{\mathrm{b}}}

\newcommand{\barcode}{\mathrm{Bar}}

\newcommand{\vect}{\mathsf{vect}}

\newcommand{\Topo}{\mathsf{Top}}

\newcommand{\Conv}{\mathrm{CH}}

\DeclareMathOperator{\id}{id}

\def\int{I}

\newcommand{\one}{\mathds{1}}

\newcommand{\Mod}{\mathrm{Mod}}
\newcommand{\Mmod}{\mathrm{mod}}

\newcommand{\Hom}{\mathrm{Hom}}

\newcommand{\Int}{\mathcal{I}}

\newcommand{\hooks}{\{\text{hooks}\}}

\newcommand{\ES}{\mathcal{E}}

\newcommand{\barcodes}{\mathrm{Barcodes}}
\newcommand{\Pers}{\mathrm{Pers}}

\newcommand{\loss}{{\mathcal L}} 

\newtheorem{theorem}{Theorem}[section]
\newtheorem{proposition}[theorem]{Proposition}

\newtheorem{corollary}[theorem]{Corollary}
\theoremstyle{definition}
\newtheorem{definition}[theorem]{Definition}
\newtheorem{remark}[theorem]{Remark}
\newtheorem{example}[theorem]{Example}

\begin{document}

\title{Differential Calculus and Optimization\\in Persistence Module Categories}

\author{Steve Oudot}
\date{}

\maketitle 
\begin{abstract}
Persistence modules are representations of products of totally ordered sets in the category
of vector spaces. They appear naturally in the representation theory of algebras, but in
recent years they have also found applications in other areas of mathematics, including
symplectic topology, complex analysis, and topological data analysis, where they arise from
filtrations of topological spaces by the sublevel sets of real-valued functions. Two
fundamental properties of persistence modules make them useful in such contexts:
(1) the fact that they are stable under perturbations of the originating functions, and 
(2) the fact that they can be approximated, in the sense of relative homological algebra, by
classes of indecomposable modules with an elementary structure.
In this text we give an introduction to the theory of persistence modules, then
we explain how the above properties can be leveraged to build a framework for differential
calculus and optimization with convergence guarantees in persistence module categories.
\end{abstract}


\section{Introduction}

Broadly speaking, persistence modules are representations of posets in the category of vector spaces over a fixed field. While these objects have been studied for decades, particularly in representation theory, the name {\em persistence modules} was coined by the topological data analysis (TDA) community at the turn of the $21^{\mathrm{st}}$ century~\cite{zomorodian2005computing}, to emphasize the role played by these modules in the study of the topology of sublevel set of functions defined on data. Meanwhile, TDA has shed a different light on these objects, raising new questions for mathematicians: a whole new theory of persistence modules  has emerged, called {\em persistence theory}~\cite{oudot2015persistence}, centered around the definition of invariants for persistence modules that can be effectively computed from data, that are provably stable under perturbations of the data, and whose constructions can be differentiated so that their parameters can be optimized in machine or deep learning contexts.

This text is a walk through persistence theory, including its basic
results and some of its latest developments. By no means is it an exhaustive account of the theory nor of its connections to other areas of mathematics: for such accounts, the reader can refer to recent textbooks on the subject~\cite{dey2022computational,polterovich2020topological,schenck2022algebraic}. The story arc here is the study of a particular class of invariants for persistence modules called {\em barcodes}, which play a central role in the $1$-parameter instance of the theory---focused on the sublevel sets of $\R$-valued functions---and which have recently been adapted as {\em signed barcodes} to the multi-parameter instance of the theory---focused on $\R^n$-valued functions. Our goal is to unveil the meaning of these invariants in terms of the structure of persistence modules, and to show how, ultimately, their construction can be differentiated and optimized in data analysis contexts.     

The text is organized as follows: Section~\ref{sec:pers_mod} introduces persistence modules and Section~\ref{sec:sublevel-sets} connects them to the topology of sublevel sets of functions; Sections~\ref{sec:structure} and~\ref{sec:stability} study their structure and stability properties, which are  leveraged in Sections~\ref{sec:diff_calc} and~\ref{sec:optim} to build a framework for differential calculus and optimization with persistence modules; the text concludes with Section~\ref{sec:perspectives} giving some directions for possible future developments on the subject.

\subsection*{Acknowledgements}
 This is a written version of the invited lecture the author gave at the $9^{\mathrm{th}}$ European Congress of Mathematics. Part of the exposition is based on the author's own work with his collaborators, whom he  wishes to thank for their invaluable contributions.

\section{Persistence Modules}
\label{sec:pers_mod}

Throughout the paper we fix a poset $(P, \leq)$  and a field $\field$. The order relation~$\leq$ is omitted in the notations whenever it is obvious from the context. We make no assumption on~$\field$, in particular it does not have to be algebraically closed and its characteristic can be arbitrary. For now, we make no assumption on~$P$ either.

\begin{definition} \label{def:pers_module}
  A {\em persistence module} over $P$ is a functor $M\colon P\to \vect_{\field}$, where $\vect_{\field}$ denotes the category of finite-dimensional vector spaces over~$\field$, and where $P$ is identified with the category having one object per element $p\in P$ and a unique morphism per pair of comparable elements $p\leq q\in\ P$. A {\em morphism} of persistence modules $M\to N$ is a natural transformation between functors $M, N\colon P\to \vect_{\field}$.
\end{definition}

Below is an example of a poset~$P$, represented by its Hasse diagram~(left),  and of a persistence module~$M$ over~$P$~(right). The diagram on the right-hand side commutes by functoriality of~$M$.

\[  
\xymatrix@=30pt{
& \bullet \\
\bullet \ar[ur] && \bullet \ar[ul]\\
&\bullet \ar[ul]\ar[ur]
}
\quad\quad\quad\quad
\xymatrix@=30pt{
& \field \ar@(ur,ul)[]|-{\, \id_{\field}\, }\\
\field^2  \ar@(dl,ul)[]|-{\text{\vbox to 10pt{\vfil\hbox to 10pt{$\id_{\field^2}$}\vfil}}} \ar[ur]|-{\text{\vbox to 10pt{\vfil\hbox to 30pt{$\left[\begin{smallmatrix}1_{\field}&1_{\field}\end{smallmatrix}\right]$}\vfil}}} && 
\field^2 \ar@(dr,ur)[]|-{\text{\vbox to 10pt{\vfil\hbox to 10pt{$\id_{\field^2}$}\vfil}}} \ar[ul]|-{\text{\vbox to 10pt{\vfil\hbox to 20pt{$\left[\begin{smallmatrix}-1_{\field}&1_{\field}\end{smallmatrix}\right]$}\vfil}}}\\
&\field  \ar@(dl,dr)[]|-{\, \id_{\field}\, } \ar[ul]|-{\text{\vbox to 15pt{\vfil\hbox to 15pt{$\left[\begin{smallmatrix}1_{\field}\\0\end{smallmatrix}\right]$}\vfil}}}\ar[ur]|-{\text{\vbox to 15pt{\vfil\hbox to 15pt{$\left[\begin{smallmatrix}0\\1_{\field}\end{smallmatrix}\right]$}\vfil}}} \ar[uu]|-{\text{\vbox to 10pt{\vfil\hbox to 10pt{$\id_{\field}$}\vfil}}}
}
\]

From now on, unless otherwise specified we assume our poset~$P$ to be a finite product of totally ordered sets, as is the case for instance in the above example:
\[ P = \prod_{i=1}^n T_i,\ \text{where each $T_i$ is totally ordered.}\]
Here $n\in\N^*$ is fixed and considered as the number of parameters of the persistence modules over~$P$, thus we talk of {\em $n$-parameter persistence modules}. Note that the totally ordered sets $T_i$ can be arbitrary, in particular they are not required to be finite.

We write $\vect_{\field}^P$ for the category of persistence modules over~$P$. As a category of functors from a small category to an abelian category,  $\vect_{\field}^P$ is itself abelian. This means in particular that it has a zero object, called the {\em zero} (or {\em trivial}) persistence module and denoted by~$0$. The category also has pointwise finite direct sums, i.e., direct sums $\bigoplus_{j\in J} M_j$ such that the set $\{ j \mid M_j(p)\neq 0\}$ is finite for every $p\in P$. Moreover, $\Hom$-sets are $\field$-vector spaces, with composition being bilinear, and every morphism has a kernel and a cokernel. The other properties of an abelian category will be used implicitly in this text.

\subsection*{Persistence modules as modules}

As their name suggests, persistence modules can be viewed as modules over a certain algebra~$\field P$ defined from the indexing poset~$P$. This algebra is freely generated by the pairs of comparable elements in~$P$:
\[ \field P \coloneqq \field^{\left(\left\{(p,q) \mid p \leq q \in P\right\}\right)}, \]
with multiplication given by:
\[ (p,q) \cdot (p',q') \coloneqq \begin{cases} (p,q') & \text{if $q=p'$}\\ 0 & \text{otherwise}, \end{cases} \]
and extended by linearity.

There is an embedding of $\vect_{\field}^P$ into the category $\Mod_{\field P}$ of right modules over~$\field P$, given by the fully faithful additive functor $\vect_{\field}^P \to \Mod_{\field P}$ defined on objects by $M\mapsto \bigoplus_{p\in P} M_p$ with right multiplication:
\[  \left(\sum_{p\in P} m_p\right)\cdot (q,r) = M(q\leq r)(m_q). \]
Thus, persistence modules can be viewed as right modules over the algebra~$\field P$. When $P$ is finite, the embedding $\vect_{\field}^P \hookrightarrow \Mod_{\field P}$ induces an equivalence of abelian categories $\vect_{\field P} \simeq \Mmod_{\field P}$, where $\Mmod_{\field P}$ is the category of finite-dimensional right $\field P$-modules.

\section{Sublevel-sets topology}
\label{sec:sublevel-sets}

  \begin{figure}[tb]
    \centering
    \includegraphics[scale=0.4]{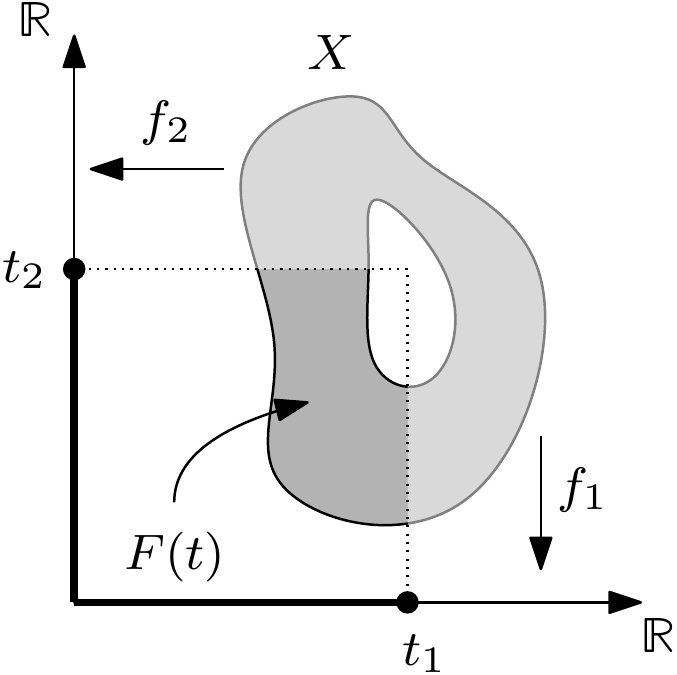}
    \caption{The $t=(t_1, t_2)$-sublevel set of the function $f=(f_1, f_2)\colon X\to\R^2$ defined by the orthogonal projections onto the coordinate axes.}
    \label{fig:sublevel-set}
  \end{figure}

  The typical use case of persistence modules is in the study of the topology of sublevel sets of functions. Given a topological space~$X$ and a function~$f=(f_1, \cdots, f_n)\colon X\to\R^n$, both of which can be arbitrary, one considers the sublevel sets of~$f$, written $F(t)$ for every $t=(t_1, \cdots, t_n)\in\R^n$ and defined as follows---see Figure~\ref{fig:sublevel-set} for an illustration:
\[ F(t) \coloneqq \bigcap_{i=1}^n f_i^{-1}\left(\left(-\infty,\,t_i\right]\right). \]
  Letting $P=\R^n$, viewed as the product of $n$ copies of the totally ordered real line and equipped with the induced product order~$\leq$, one gets:
  \[ F(s) \subseteq F(t)\quad \forall s\leq t \in\R^n. \] 
  This yields a functor $F\colon P\to\Topo$, where $\Topo$ denotes the category of topological spaces, defined as above on objects and sending each morphism $s\leq t\in \R^n$ to the inclusion map $F(s)\hookrightarrow F(t)$. This functor is called the {\em filtration} of the sublevel sets of~$f$. Assuming (as we do from now on) that the sublevel sets of~$f$ each have finite-dimensional singular homology in every degree, one post-composes~$F$ by the singular homology functor~$H_*$ in some fixed degree~$*$ with cofficients in~$\field$, to get a persistence module $H_*(f)\colon \R^n\to\vect_{\field}$ called the {\em persistent homology of~$f$ in degree~$*$}.

  Let us give two classic examples in the context of TDA.

\begin{example} \label{ex:filt_distance}
  Let $A\subset\R^d$ be a finite set of data points. Take $X=\R^d$ and define $f\colon \R^d\to\R$ as the distance to the data:
  \[ f(x) \coloneqq \min_{a\in A} \left\|x-a\right\|_2. \]
Its sublevel sets are unions of Euclidean balls of same radius around the data points, as illustrated in Figure~\ref{fig:filt_dist}.
  \begin{figure}[tb]
    \centering
    \includegraphics[width=1\textwidth]{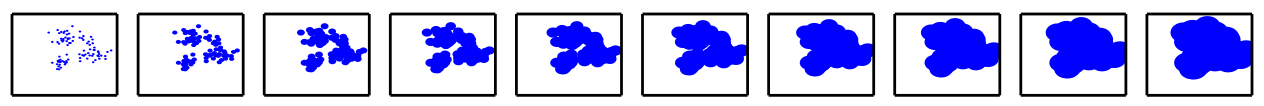}
    \caption{Filtration of the distance function (excerpt with $10$ different levels~$t$, starting at $t=0$).}
    \label{fig:filt_dist}
    \end{figure}
  The corresponding $1$-parameter persistence modules in singular homology of degrees~$0$ and~$1$ respectively over~$\field$, noted $H_0(f)$ and $H_1(f)$ respectively, are shown below (restricted to the same $10$ indices~$t$ as in Figure~\ref{fig:filt_dist}):
  \begin{align*}
    H_0(f)\colon & \xymatrix@C=20pt{
      \field^{84} \ar@{->>}[r] & \field^9 \ar@{->>}[r] & \field^3 \ar@{->>}[r] & \field^2 \ar@{->>}[r] & \field \ar[r]^-{\simeq} & \field \ar[r]^-{\simeq} & \field \ar[r]^-{\simeq} & \field \ar[r]^-{\simeq} & \field \ar[r]^-{\simeq} & \field } \\
    H_1(f)\colon & \xymatrix@C=22pt{
0 \ar[r] & \field^3 \ar[r]^-0 & \field \ar[r]^-0 & \field \ar[r]^-0 & \field^2 \ar[r]^-{\left[\begin{smallmatrix}1&0\end{smallmatrix}\right]} & \field \ar[r] & 0 \ar[r] &0 \ar[r] &0 \ar[r] &0 }
  \end{align*}  
\end{example}

\begin{example} \label{ex:filt_distance-density}
  Let again $A\subset\R^d$ be a finite set of data points. Take $X=\R^d$ and define $f=(f_1, f_2)\colon \R^d\to\R^2$ as the combination of the distance to the data with some co-density estimator:
  \begin{align*} f_1(x) & \coloneqq \min_{a\in A} \left\|x-a\right\|_2, \\
    f_2(x) & \coloneqq \#\left\{ a\in A \mid \left\|x-a\right\|_2 >\e \right\} \quad \text{for some fixed parameter~$\e\geq 0$}.
  \end{align*}
Its sublevel sets are unions of Euclidean balls of same radius around data points with high enough density, as illustrated in Figure~\ref{fig:filt_dist_dens}.
  \begin{figure}[tb]
    \centering
    \includegraphics[width=1\textwidth]{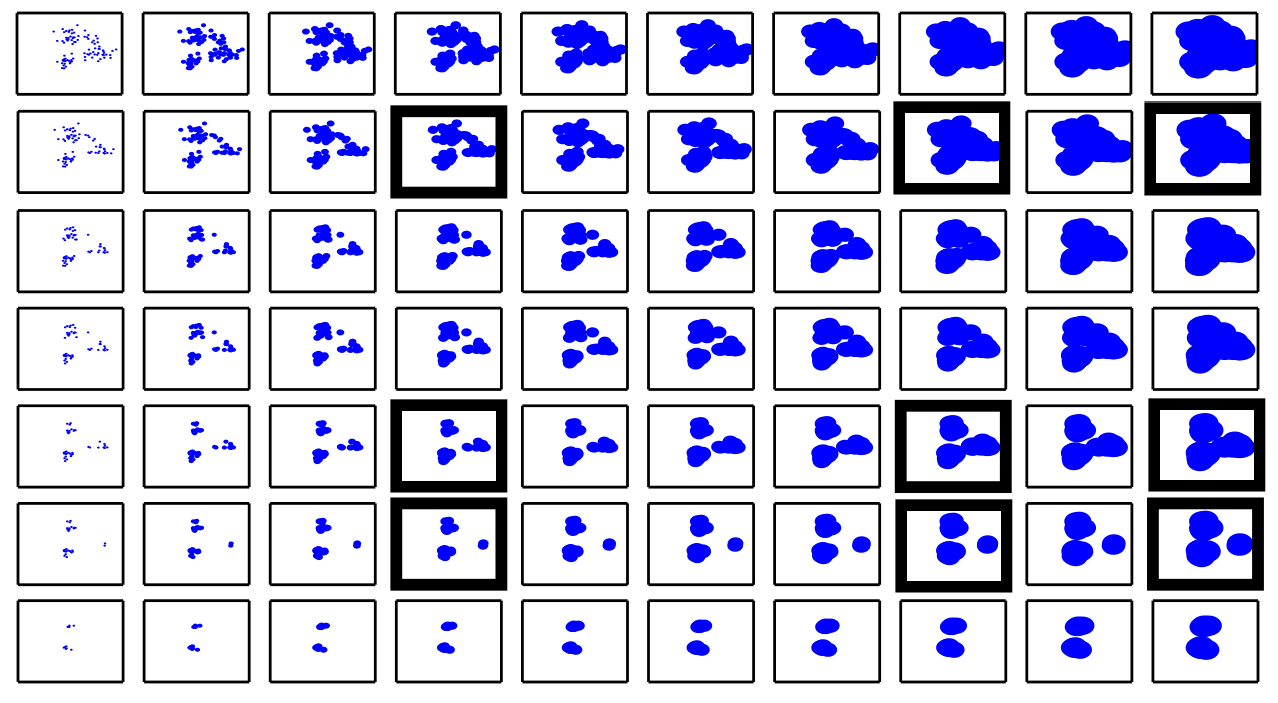}
    \caption{Bi-filtration of the distance function and co-density estimator (excerpt with $10$ different levels for the distance and $7$ for the co-density). The distance levels are represented horizontally, the co-density levels vertically. The indices highlighted by thick rectangles arranged in a $3\times 3$ grid exhibit three connected components that merge two-by-two in three different ways then all together.}
    \label{fig:filt_dist_dens}
    \end{figure}
  The corresponding $2$-parameter persistence module in singular homology of degree~$0$ over~$\field$ restricts to the following persistence module over the $3\times 3$ grid highlighted in Figure~\ref{fig:filt_dist_dens}:
  \begin{equation}\label{eq:indec_3x3}
    \begin{gathered}
    \xymatrix@=25pt{
  \field^2 \ar^-{\left[\begin{smallmatrix}1&1\end{smallmatrix}\right]}[r] & \field \ar^-{\id}[r] & \field\\
  \field^3 \ar^-{\left[\begin{smallmatrix}1&0&0\\0&1&1\end{smallmatrix}\right]}[u]
  \ar^-{\left[\begin{smallmatrix}1&0&1\\0&1&0\end{smallmatrix}\right]}[r] &
  \field^2
  \ar_-{\left[\begin{smallmatrix}1&1\end{smallmatrix}\right]}[u]
  \ar^-{\left[\begin{smallmatrix}1&1\end{smallmatrix}\right]}[r] &
  \field \ar_-{\id}[u]\\
  \field^3 \ar^-{\id}[u] \ar_-{\id}[r] & \field^3
  \ar_-{\left[\begin{smallmatrix}1&0&1\\0&1&0\end{smallmatrix}\right]}[u]
  \ar_-{\left[\begin{smallmatrix}1&1&0\\0&0&1\end{smallmatrix}\right]}[r] & \field^2 \ar_-{\left[\begin{smallmatrix}1&1\end{smallmatrix}\right]}[u]
    }
    \end{gathered}
    \end{equation}
\end{example}

  \begin{figure}[t]
    \centering
    \includegraphics[width=1\textwidth]{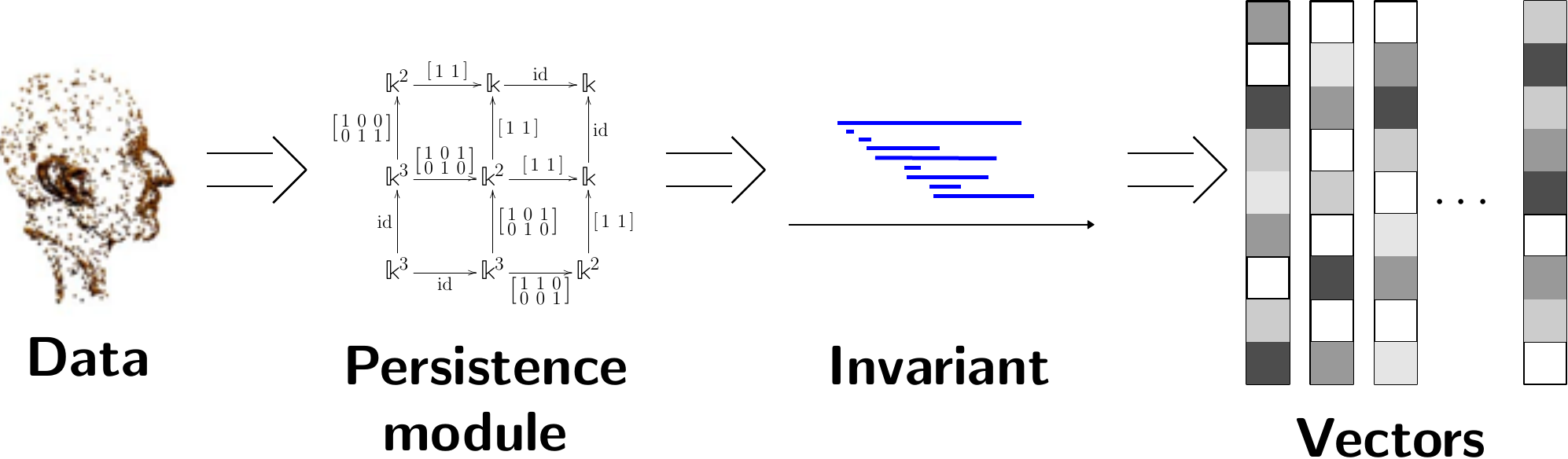}
    \caption{The topological data analysis pipeline.}
    \label{fig:TDA_pipeline}
    \end{figure}

  More broadly, in the context of TDA, one builds persistence modules from data, similarly to Examples~\ref{ex:filt_distance} and~\ref{ex:filt_distance-density}. Then, one computes invariants on these modules: a variety of invariants have been proposed, in this text they all take the form of {\em barcodes}, i.e., multisets of intervals in the indexing poset. Then, one turns these barcodes into vectors in order to feed them to some machine of deep learning architecture. The corresponding pipeline is illustrated in Figure~\ref{fig:TDA_pipeline}. Two fundamental properties of persistence modules make this pipeline relevant for data analysis:
  \begin{itemize}
    \setlength{\itemindent}{10pt}
  \item[\textsc{(structure)}] persistence modules can be decomposed, or approximated in the sense of homological algebra, by certain classes of modules with an elementary structure called {\em interval modules}; this is where their barcodes are defined  from;
    \item[\textsc{(stability)}] persistence modules and their barcodes are stable with respect to perturbations of their originating functions; this enables the use of barcodes as stable features for data.
  \end{itemize}
  The next two sections are devoted to fleshing out these properties:  Section~\ref{sec:structure} for the structure, Section~\ref{sec:stability} for the stability.

  \section{Structure of persistence modules}
  \label{sec:structure}

  As before, let $(P,\leq)$ be a poset and $\field$ a field. Unless otherwise stated, we make no assumption on~$P$. There are a variety of ways to encode the structure of persistence modules over~$P$, here we report on two of them: via direct-sum decompositions (Section~\ref{sec:direct-sum}), and via homological approximations (Section~\ref{sec:homological-approx}). In both cases, interval modules play a central part in the story, so we begin by giving their definition.

  \begin{definition}\label{def:int}
    An {\em interval of $P$} is a subset $I\subseteq P$ that is both convex and connected with respect to the order~$\leq$ on~$P$, that is:
  \begin{align*}
    \text{\textsc{(convex)}} & \quad s,t\in  I\ \implies\ u\in I\quad \forall\ s\leq u\leq t, \\ 
    \text{\textsc{(connected)}} & \quad s,t\in  I\ \implies\ \exists r\in\N^*,\ \exists \{u_i\}_{i=0}^r \subseteq I\ \colon \ s=u_0\leq u_1\geq \cdots \geq  u_r=t.
  \end{align*}
  \end{definition}

  \begin{definition}\label{def:int_mod}
    Given an interval $I\subseteq P$, the corresponding {\em interval module ${\field}_\int$} is defined by:
      \[
        \field_I(t) = \left\{ 
\begin{array}{l}
\field \ \mbox{if}\ t\in I\\[0.5em]
0\ \mbox{otherwise}
\end{array}
\right.
\quad\quad\quad
\field_I(s\leq t) = \left\{ 
\begin{array}{l}
\id_{\field} \ \mbox{if}\ s,t\in I\\[0.5em]
0\ \mbox{otherwise}
\end{array}
\right.
 \]
\end{definition}

  Interval modules are useful for applications because they are entirely characterized by their support, which is an interval of~$P$ and therefore a geometric object.  This makes it possible for practitioners to think purely in terms of the geometry of these intervals, thus forgetting about the algebra. Examples will be given in the following.
  
  \subsection{Direct-sum decompositions}
  \label{sec:direct-sum}

A persistence module~$M$ is said to be {\em decomposable} if it decomposes as $M\simeq N \oplus N'$ for some modules $N, N' \neq 0$. Otherwise, if no such decomposition exists, $M$ is called {\em indecomposable}. 
Note that interval modules  are indecomposable: this is because their endomorphism ring is isomorphic to~$\field$, hence local. Therefore, interval modules can serve as building blocks for decomposing arbitrary persistence modules. The hope is that no other types of indecomposables appear in the decomposition, which is  what happens for instance when the poset $P$ is totally ordered, an important special case that encompasses the study of the topology of the sublevel sets of real-valued functions as in Example~\ref{ex:filt_distance}.  

  \begin{theorem}[\cite{botnan-crawleybovey,Crawley-Boevey2012}]
    \label{th:decomp_total-order}
    Suppose $(P, \leq)$ is totally ordered. Then, any $M\in \vect_{\field}^P$ decomposes as follows:
    \[ M \simeq \bigoplus_{j\in J}\, \field_{I_j}, \]
    where each $I_j$ is an interval of~$P$. This decomposition is unique up to isomorphism and reordering of the terms.
  \end{theorem}
  
When $P$ is finite, this is just Gabriel's theorem~\cite{Gabriel1972} specialized to the case of $A_n$ type quivers.  The version of the result stated here, due to Botnan and Crawley-Boevey, makes no assumption on~$P$ (beside being totally ordered) nor on~$M$.

\begin{example}\label{ex:filt_distance_decomp}
  The persistence modules from Example~\ref{ex:filt_distance} decompose as follows:
  \begin{align*}
    H_0(f) & \simeq  \xymatrix@C=20pt{
      \field \ar^-{\id_{\field}}[r] & \field \ar^-{\id_{\field}}[r] & \field \ar^-{\id_{\field}}[r] & \field \ar^-{\id_{\field}}[r] & \field \ar[r]^-{\id_{\field}} & \field \ar[r]^-{\id_{\field}} & \field \ar[r]^-{\id_{\field}} & \field \ar[r]^-{\id_{\field}} & \field \ar[r]^-{\id_{\field}} & \field }\\
    & \oplus \xymatrix@C=20pt{
      \field \ar^-{\id_{\field}}[r] & \field \ar^-{\id_{\field}}[r] & \field \ar^-{\id_{\field}}[r] & \field \ar[r] &  0 \ar[r] &  0 \ar[r] &  0 \ar[r] &  0 \ar[r] &  0 \ar[r] &  0 }\\
    & \oplus \xymatrix@C=20pt{
      \field \ar^-{\id_{\field}}[r] & \field \ar^-{\id_{\field}}[r] & \field \ar[r] & 0 \ar[r] &  0 \ar[r] &  0 \ar[r] &  0 \ar[r] &  0 \ar[r] &  0 \ar[r] &  0}\\
    & \oplus \left(\xymatrix@C=20pt{
      \field \ar^-{\id_{\field}}[r] & \field \ar[r] & 0 \ar[r] & 0 \ar[r] &  0 \ar[r] &  0 \ar[r] &  0 \ar[r] &  0 \ar[r] &  0 \ar[r] &  0}\right)^6\\
    & \oplus \left(\xymatrix@C=20pt{
      \field \ar[r] & 0 \ar[r] & 0 \ar[r] & 0 \ar[r] &  0 \ar[r] &  0 \ar[r] &  0 \ar[r] &  0 \ar[r] &  0 \ar[r] &  0}\right)^{75};\\[1.5ex]
H_1(f) & \simeq
    \xymatrix@C=20pt{
      0 \ar[r] & 0 \ar[r] & 0 \ar[r] & 0 \ar[r] & \field \ar^-{\id_{\field}}[r] & \field \ar[r] & 0 \ar[r] & 0 \ar[r] &  0 \ar[r] &  0} \\
    & \oplus \xymatrix@C=20pt{
     0 \ar[r] & 0 \ar[r] & 0 \ar[r] & 0 \ar[r] & \field \ar[r] & 0 \ar[r] & 0 \ar[r] & 0 \ar[r] &  0 \ar[r] &  0}\\
    & \oplus \xymatrix@C=20pt{
     0 \ar[r] & 0 \ar[r] & 0 \ar[r] & \field \ar[r] & 0 \ar[r] & 0 \ar[r] & 0 \ar[r] &  0 \ar[r] &  0 \ar[r] &  0}\\
    & \oplus \xymatrix@C=20pt{
      0 \ar[r] & 0 \ar[r] & \field \ar[r] & 0 \ar[r] & 0 \ar[r] & 0 \ar[r] &  0 \ar[r] &  0 \ar[r] &  0 \ar[r] &  0} \\
   & \left( \xymatrix@C=20pt{
     0 \ar[r] & \field \ar[r] & 0 \ar[r] & 0 \ar[r] & 0 \ar[r] &  0 \ar[r] &  0 \ar[r] &  0 \ar[r] &  0 \ar[r] &  0} \right)^3.
  \end{align*}
\end{example}

Interval decompositions, when they exist,  provide a direct interpretation of the structure of the persistence modules: in Example~\ref{ex:filt_distance_decomp} for instance, the supports of the interval summands of~$H_0(f)$ encode the lifespan of each path-connected component in the union of balls as the radius grows; meanwhile, all but one of the interval summands of~$H_1(f)$ are supported on singletons, indicating that the corresponding holes appearing in the growing union of balls have an ephemeral existence, while the non-singleton interval has a corresponding hole (the main one in the picture of Figure~\ref{fig:filt_dist}) that lasts over two levels in the excerpt. This motivates the use of interval decompositions---when they exist---in applications. Or, rather, the use of their associated {\em barcodes}, which encode the (unique) geometric information contained in each decomposition, that is, the supports of the interval summands.

\begin{definition}\label{def:barcode}
  Suppose a persistence module $M$ over an arbitrary poset~$P$
  decomposes as a direct sum of interval modules:
  \[ M \simeq \bigoplus_{j\in J} \field_{I_j},\ \text{where each $I_j$ is an interval of~$P$}. \]
  Then, $M$ is called {\em interval decomposable}, and by uniqueness of the decomposition up to isomorphism and reordering of the terms, the multiset of intervals involved in the decomposition is unique. This multiset is called the {\em barcode} of~$M$, noted $\barcode(M)$:
  \[ \barcode(M) \coloneqq \left\{ I_j \mid j \in J\right\}. \]
\end{definition}

\begin{example}
  The interval decompositions of Example~\ref{ex:filt_distance_decomp} lead to the following barcodes (in blue):
  \[
  \xymatrix@R=2pt{\\\\\barcode (H_0(f)) =\\\\}
  \xymatrix@R=2pt@C=20pt{
\bullet \ar[r] & \bullet \ar[r] & \bullet \ar[r] & \bullet \ar[r] & \bullet \ar[r] & \bullet \ar[r] & \bullet \ar[r] & \bullet \ar[r] & \bullet \ar[r] & \bullet \\
\ar@{{*}-{*}}@[blue][rrrrrrrrr] &&&&&&&&& \\
\ar@{{*}-{*}}@[blue][rrr] &&&\\
\ar@{{*}-{*}}@[blue][rr] &&\\
\ar_>{\color{blue}_{(6)}}@{{*}-{*}}@[blue][r] & \\
\ar_<{{\color{blue}_{(75)}}}@{{*}-}@[blue][r]|{~\hspace{30pt}} &
  }\]
  \[
  \xymatrix@R=2pt{\\\\\barcode (H_1(f)) =\\\\}
  \xymatrix@R=2pt@C=20pt{
\bullet \ar[r] & \bullet \ar[r] & \bullet \ar[r] & \bullet \ar[r] & \bullet \ar[r] & \bullet \ar[r] & \bullet \ar[r] & \bullet \ar[r] & \bullet \ar[r] & \bullet \\
&&&& \ar@{{*}-{*}}@[blue][r] & \\
&&&& \ar@{{*}-}@[blue][r]|{~\hspace{30pt}} &\\
&&& \ar@{{*}-}@[blue][r]|{~\hspace{30pt}} &\\
&& \ar@{{*}-}@[blue][r]|{~\hspace{30pt}} &\\
& \ar_<{{\color{blue}_{(3)}}}@{{*}-}@[blue][r]|{~\hspace{30pt}} &
  }\]
\end{example}

Unfortunately, when $P$ is not totally ordered---for instance when  $P$ is a finite product of $n>1$ totally ordered sets, interval decomposability is no longer guaranteed, even though every persistence module~$M$ still decomposes into  indecomposable summands.
  \begin{theorem}[\cite{botnan-crawleybovey}]
    \label{th:decomp}
    For any poset $(P, \leq)$, any $M\in \vect_{\field}^P$ decomposes as follows:
    \[ M \simeq \bigoplus_{j\in J}\, M_j, \]
    where each summand $M_j$ is indecomposable with local endomorphism ring. This decomposition is unique up to isomorphism and reordering of the terms.
  \end{theorem}
  The summands~$M_j$ appearing in the decomposition may not be interval modules: for instance, the persistence module~\eqref{eq:indec_3x3} is indecomposable yet not an interval module, being not of pointwise dimension~$\leq 1$. And such summands do appear in applications, since~\eqref{eq:indec_3x3} is coming from Example~\ref{ex:filt_distance-density}. From the algebraic standpoint, the situation is even worse: not only may persistence modules have non-interval summands, but in fact the algebra~$\field P$ is generally of {\em wild representation type}, since for instance one can embed the representation category of the following wild-type quiver into the representation category of $\field P$ when $P$ contains a large enough $n$-dimensional ($n>1$) grid:
  \[ \xymatrix@=10pt{
    \bullet \ar[rrrr] &&&& \bullet\\
    & \bullet \ar[urrr]\\
    && \bullet \ar[uurr]\\
    &&& \bullet \ar[uuur]\\
    &&&& \bullet \ar[uuuu]
    }\]
For us, this means that, unless $P$ is very small (not more than a handful of elements), when $n>1$ there is no canonical classification of all the indecomposable persistence modules over the poset~$P$, therefore parametrizing the indecomposables to build up unique barcodes from direct-sum decompositions is an impossible task.

\subsection{Homological approximations}
  \label{sec:homological-approx}

  Given the above shortcoming, and the importance of working with invariants that can be interpreted in application contexts, the TDA community has made significant effort in the past decade to propose alternative invariants from which meaningful barcodes can be built. These invariants are necessarily incomplete, considering the wild representation type of~$\field P$. Here we describe one class of such invariants, based on homological approximations of the modules.

  From now on, we assume that the poset $P$ is finite, so its algebra~$\field P$ is finite-dimensional. We can then rely on the well-established theory of representations of Artin algebras. The indecomposable projectives in this setting are the interval modules $\field_{[p, \infty)}$, where $[p, \infty)\coloneqq  \left\{q\in P \mid q\geq p\right\}$ is the principal upset at~$p$ in~$P$. Thus, by taking successive projective approximations, we can build left resolutions of any persistence module~$M$ by interval modules of this type. Moreover, assuming~$M$ admits finite projective presentations, the minimal such resolution has finite length, bounded by the cardinality of~$P$---or, more tightly, by the number~$n$ of its parameters when $P$ is a product of $n$ finite totally ordered sets, by Hilbert's Syzygy theorem:
      \[\xymatrix{
0\ar[r]& 
\displaystyle\bigoplus_{p\in J_n} \field_{[p, \infty)} \ar[r] & 
\displaystyle\bigoplus_{p\in J_{n-1}} \field_{[p, \infty)} \ar[r] & \cdots\ar[r] &
  \displaystyle\bigoplus_{p\in J_0} \field_{[p, \infty)} \ar@{->>}[r] & M}\]
As the minimal resolution is unique up to isomorphism, the supports  of the summands in the various terms of the resolution are uniquely defined. We can then gather them into a barcode: more precisely, for reasons that we will explain afterwards, we build a composite object called the {\em signed barcode} of~$M$, noted $\barcode_{\llcorner}(M)$,  composed of two regular barcodes, the first one (called the {\em positive part} and noted $\barcode_{\llcorner}^+(M)$) gathers the supports coming from even degrees in the resolution, and the second one (called the {\em negative part} and noted $\barcode_{\llcorner}^-(M)$)  gathers the supports coming from odd degrees:
\begin{equation}\label{eq:signed_barcode_usual_proj}
  \begin{array}{rl}
    \barcode_{\llcorner}(M) \coloneqq & \left(\barcode_{\llcorner}^+(M),\, \barcode_{\llcorner}^-(M)\right)\\[1em]
    \text{where}&\begin{cases}\barcode_{\llcorner}^+(M)\coloneqq\bigsqcup_{i\in 2\N}\, \left\{[p, \infty) \mid p\in J_i\right\}\\\barcode_{\llcorner}^-(M)\coloneqq\bigsqcup_{i\in 2\N+1}\, \left\{[p, \infty) \mid p\in J_i\right\}.\end{cases}
\end{array}
\end{equation}
Compared to regular barcodes, signed barcodes put signs on their intervals based on the parity of the degree in the resolution in which they appear. This twist is motivated by what happens in the Grothendieck group~$K_0$, where the resolution yields the following signed decomposition of the equivalence class~$\left[M\right]$ of~$M$:
\[ \left[M\right]\ =\ 
\sum\limits_{i\in \N}\,  (-1)^i \left(\sum_{p\in J_i} \left[{\field}_{[p, \infty)}\right]\right) 
\ =\ 
\sum\limits_{p\in \bigsqcup\limits_{i\ \text{even}} J_i} \left[{\field}_{[p, \infty)}\right] 
- \sum\limits_{p\in \bigsqcup\limits_{i\ \text{odd}} J_i} \left[{\field}_{[p, \infty)}\right].
\]
Note that some terms may cancel out in this alternating sum, in which case they disappear from the decomposition of~$\left[M\right]$ in~$K_0$. By contrast, we keep them  in the positive and negative parts of the signed barcode~\eqref{eq:signed_barcode_usual_proj}, thus preserving all the information contained in the terms of the minimal projective resolution of~$M$. In practice this is more costly to compute, since we have to go through computing the resolution as opposed to directly decomposing $M$ in $K_0$. But it yields a stronger invariant, with provable stability properties as we shall see in the next section.

The biggest limitation of this approach is that it does not make use of the full set of intervals in~$P$ to approximate the modules,  restricting itself to using the subset of principal upsets. In order to enable the use of larger collections of intervals to build our resolutions, we rely on the theory of relative homological algebra for Artin algebras developed in~\cite{auslander-solberg}. Let thus $\Int$ be a set of intervals in~$P$. In view of the above discussion, we assume that $\Int$ contains all the principal upsets $[p, \infty)$ for $p\in P$. The bottomline is to make the interval modules $\field_{\int}$ for $\int\in\Int$ projective. For this, we consider only those short exact sequences on which all the $\Hom$ functors $\Hom(\field_\int, -)$ for $\int\in\Int$ are exact.
  Let $\ES_{\Int}$ denote the collection of such sequences:
\begin{align*}
  \ES_{\Int} \coloneqq & \left\{ 0\to L\to M \to N \to 0 \ \text{exact}\ \mid \right.\\
  & \left. \ 0\to \Hom(\field_{\int}, L)\to \Hom(\field_{\int}, M) \to \Hom(\field_{\int}, N) \to 0 \ \text{exact for all}\ \int\in\Int\right\}.
\end{align*}

\begin{proposition}[\cite{auslander-solberg}]
  \label{prop:exact_struct}
    ~
\begin{itemize}
\item[(i)] $\vect_{\field}^{P}$ equipped with $\ES_\Int$ has the structure of an exact category;
\item[(ii)] $\vect_{\field}^{P}$ has enough projectives relative to $\ES_\Int$; 
\item [(iii)] the projectives relative to $\ES_\Int$ are the finite direct sums of interval modules of type $\field_\int$ for $\int\in\Int$.
\end{itemize}
\end{proposition}

Item~(i) implies that we may consider projective resolutions of persistence modules relative to~$\ES_\Int$, when they exist. Item~(ii) comes from our assumption that $\Int$ contains all the principal upsets in~$P$, so the projectives relative to~$\ES_\Int$ include in particular the usual projectives. This item ensures that every persistence module does admit projective resolutions relative to~$\ES_\Int$. Finally, item~(iii) guarantees that the terms in such resolutions are finite direct sums of interval modules supported on elements of~$\Int$ as desired.  
We can then take projective resolutions of any persistence module~$M$ relative to~$\ES_\Int$. More precisely, we take the minimal such resolution:
\[ \xymatrix{
0 \ar[r] & \displaystyle\bigoplus_{I\in J_r} \field_{I} \ar[r] & \displaystyle\bigoplus_{I\in J_{r-1}} \field_{I} \ar[r] & \cdots \ar[r] &
 \displaystyle\bigoplus_{I\in J_1} \field_{I} \ar[r] & \displaystyle\bigoplus_{I\in J_0} \field_{I} \ar@{->>}[r] & M} \]
Provided this resolution has finite length, we can form a signed barcode as before, but this time relative to our choice of family~$\Int$ of intervals:
\begin{equation}\label{eq:signed_barcode_rel_proj}
  \barcode_{\Int}(M) \coloneqq 
  \left(\barcode_{\Int}^+(M),\, \barcode_{\Int}^-(M)\right)
  \quad \text{where} \quad 
  \begin{cases}\barcode_{\Int}^+(M)\coloneqq\bigsqcup_{i\in 2\N}\, J_i\\\barcode_{\Int}^-(M)\coloneqq\bigsqcup_{i\in 2\N+1}\, J_i.\end{cases}
\end{equation}
Again, the signs in the barcode match with the ones in the decomposition of the equivalence class $\left[ M\right]$ of~$M$  in the Grothendieck group $K_0^{\ES_{\Int}}$ relative to $\ES_{\int}$, and the difference is that the signed barcode keeps the terms that cancel out in~$K_0^{\ES_{\Int}}$:
\begin{equation}\label{eq:dec_Grothendieck}
  \ \left[M\right]\ =\ 
\sum\limits_{i\in \N}\,  (-1)^i \left(\sum_{\int\in J_i} \left[{\field}_{\int}\right]\right) 
\ =\ 
\sum\limits_{\int\in \bigsqcup\limits_{i\ \text{even}} J_i} \left[{\field}_{\int}\right] 
- \sum\limits_{\int\in \bigsqcup\limits_{i\ \text{odd}} J_i} \left[{\field}_{\int}\right]. 
\end{equation}

\begin{example}\label{ex:hook-resol}
  Consider the indecomposable module~$M$ from~\eqref{eq:indec_3x3}. Let $\Int = \{ [p,\infty) \setminus [q,\infty) \mid p<q\in P \cup \{\infty\}\}$. The minimal projective resolution of~$M$ relative to~$\ES_{\Int}$,which has length~$1$, is shown below without specifying the differentials (which are irrelevant to our invariant) and with the supports of the summands in the resolution colored in blue or in red depending on the parity of their degree in the resolution:
      \begin{center}
        \includegraphics[width=\textwidth]{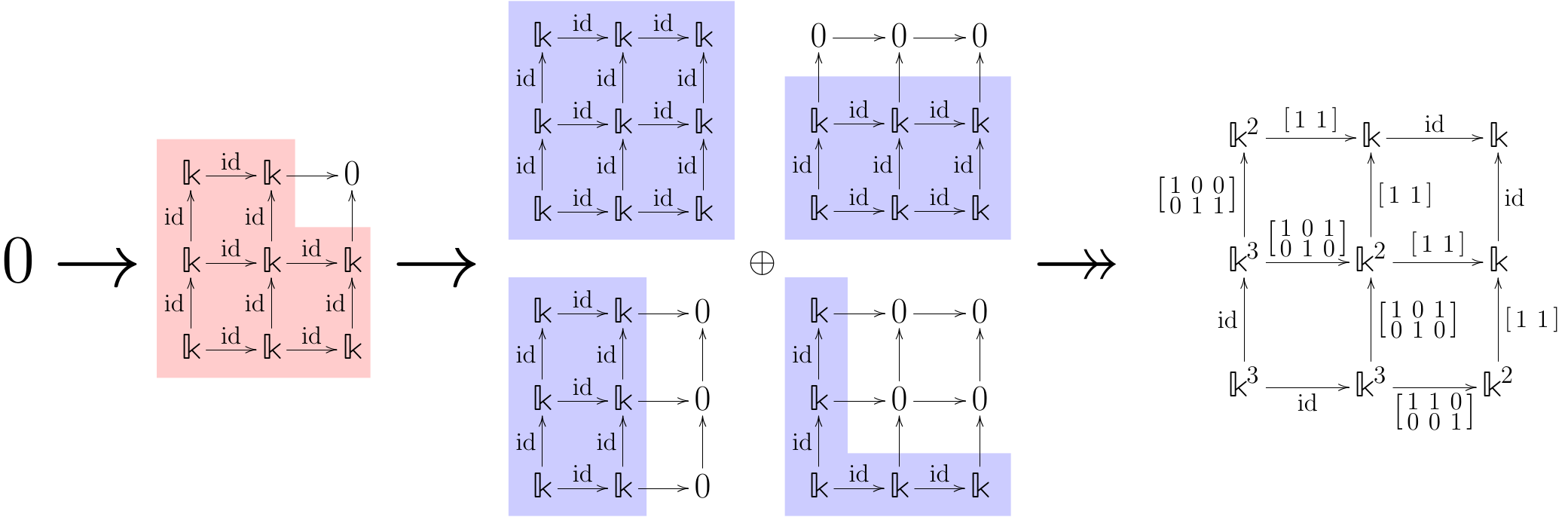}
      \end{center}
      The blue supports form the positive part of the signed barcode $\barcode_{\Int}(M)$, while the red supports form the negative part---see Figure~\ref{fig:indec_signed-barcode} for an alternative representation of $\barcode_{\Int}(M)$ as an actual barcode in the Euclidean plane.
      \begin{figure}[tb]
        \centering
        \includegraphics[scale=0.4]{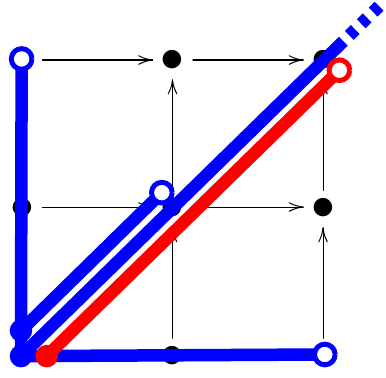}
        \caption{The signed barcode~$\barcode_{\Int}(M)$ of the indecomposable module~$M$ of Example~\ref{ex:hook-resol}, represented as an actual barcode in the Euclidean plane, where, by convention, each interval $[p,\infty) \setminus [q,\infty)$ gives rise to a copy of the line segment $[p,q)$, while each interval $[p,\infty)$ gives rise to a copy of the diagonal ray starting at~$p$. The bars in the barcode are colored according to their sign: blue for positive, red for negative.}
        \label{fig:indec_signed-barcode}
      \end{figure}
      Note that, if we replaced $M$ by  $M\oplus \field_{\int_1}$, where $\field_{\int_1}$ is the interval module appearing in degree~$1$ in the above resolution, then the minimal projective resolution would gain one copy of~$\field_{\int_1}$ in degree~$0$, which would be added to the positive part of the signed barcode but would cancel out with the copy of~$\field_{\int_1}$ in degree~$1$ in the relative Grothendieck group~$K_0^{\ES_{\Int}}$.
\end{example}

The intervals considered in the above example are called {\em hooks}, because when $P$ has $2$ parameters they have the shapes of (possibly degenerate) hooks, as we saw in the example. We therefore write $\Int =\hooks$ for their family. These intervals are of particular interest to us because their corresponding interval modules are quotients of two usual indecomposable projectives, so in this sense $\hooks$ is the simplest family of intervals to consider beyond the principal upsets. Another related consideration is that the intervals that appear in the direct-sum decomposition of finitely presented modules in the 1-parameter setting (Theorem~\ref{th:decomp_total-order}) are right-open, therefore  they are precisely the 1-dimensional hooks.   

For the previous construction to make sense, it is important that the minimal projective resolution of~$M$ relative to~$\ES_\Int$ has finite length, for otherwise the alternating sum in~\eqref{eq:dec_Grothendieck} is not defined. Another motivation is that, in practice, one  can only store finite barcodes on a computer. This is why bounding the global dimension  $\mathrm{gldim}^{\ES_\Int}\!\!\left(\vect_{\field}^{P}\right)$ of the persistence module category relative to the chosen set~$\Int$ of intervals has become a prominent question for the TDA community in recent years. Here are some of the results obtained:

\begin{theorem}[\cite{asashiba-escolar-nakashima-yoshiwaki,blanchette-brustle-hanson,botnan-oppermann-oudot,botnan2024bottleneck,cacholski-guidolin-ren-scolamiero-tombari,miller}]
  \label{th:bounds_gldim}
  ~\\
  For any finite poset~$(P,\leq)$, one has $\mathrm{gldim}^{\ES_\Int}\!\!\left(\vect_{\field}^{P}\right)<\infty$ for the following sets~$\Int$ of intervals:
  \begin{itemize}
\item[$\bullet$] $\Int = \hooks$, 
\item[$\bullet$] $\Int = \{\text{segments}\} \coloneqq \{ [p,q] \mid p\leq q\in P\}$,
\item[$\bullet$] $\Int = \{\text{upsets}\} \coloneqq \{ [Q, \infty) \mid Q\subseteq P\}$, where  $[Q, \infty) \coloneqq \{ s\in P \mid \exists q\in Q\ \text{s.t.}\ s\geq q\}$,
\item[$\bullet$] $\Int = \{\text{single-source intervals}\} \coloneqq \{ [p, \infty) \cap (-\infty, Q] \mid p\in P,\, Q\subseteq P\}$, where $(-\infty, Q] \coloneqq \{ s\in P \mid \exists q\in Q\ \text{s.t.}\ s\leq q\}$,
\item[$\bullet$] $\Int = \{ \text{intervals}\}$.
\end{itemize}
\end{theorem}

Almost all the upper bounds on the global dimension obtained so far depend on the size of the poset~$P$, therefore they do not yield any interesting upper bound when $P$ is infinite, for instance when $P=\R^n$. A notable exception is the bound obtained for hooks, which depends only on the number~$n$ of parameters of~$P$. More precisely, it was proven  in~\cite{botnan2024bottleneck}  that $\mathrm{gldim}^{\ES_{\hooks}}\!\!\left(\vect_{\field}^{P}\right)$ is exactly~$2n-2$ when $P=\prod_{i=1}^n T_i$ with each $T_i$ finite and totally ordered. This linear dependence in~$n$ of the global dimension is remarkable in that it is of the same order as its counterpart in the usual projective setting, where  Hilbert's Syzygy theorem implies that $\mathrm{gldim}\!\left(\vect_{\field}^{P}\right) = n$.  The result was then naturally extended to the case $P=\R^n$, assuming the modules under consideration are finitely presented---and denoting by~$\vect_{\field, \mathrm{fp}}^{\R^n}$ the corresponding full subcategory:
\begin{theorem}[\cite{botnan2024bottleneck}]
  \label{th:gldim_Rn}
 $\mathrm{gldim}^{\ES_{\hooks}}\!\!\left(\vect_{\field, \mathrm{fp}}^{\R^n}\right) = 2n-2$.
\end{theorem}
Note that this result recovers the interval decomposition theorem~\ref{th:decomp_total-order} for finitely presented 1-parameter persistence modules as a special case: ineed, when $n=1$ the global dimension relative to hooks is~$2-2=0$, which means that every finitely presented persistence module over~$\R$ is projective relative to~$\ES_{\hooks}$ hence interval decomposable with right-open (i.e. hook) intervals. A more general and systematic treatment of the connection between the finite and infinite poset settings can be found in~\cite{amiot2024invariants,blanchette-brustle-hanson-2}.

  \section{Stability of persistence modules}
  \label{sec:stability}

Throughout this section we assume the ground poset~$P$ to be~$\R^n$ equipped with the product order---noted~$\leq$. 
  
  \subsection{Interleaving distance}
  \label{sec:interleaving}

  Expressing the stability of persistence modules requires first to equip their category $\vect_{\field}^{\R^n}$ with a metric. Among the many choices of metrics available, one turns out to be universal  in the context of sublevel-sets topology: the {\em interleaving distance}. In some sense, this distance measures the defect of isomorphism between persistence modules. In order to define it, we first introduce the  {\em $\e$-shift functor} $-[\e]$, parametrized by $\e\geq 0$ and defined on objects by:
  \begin{align*} M[\e](t)& \coloneqq M\left(t+\e\one\right), \\
    M[\e](s\leq t)& \coloneqq M\left(s+\e\one \leq t+\e\one\right),\\
    \text{where}\ \one & \coloneqq \left[\begin{smallmatrix}1\\\vdots\\[0.8ex]1\end{smallmatrix}\right]\in\R^n,
  \end{align*}
  and likewise on morphisms. For any object~$M$, there is a natural morphism $\eta^\e_M\colon M\to M[\e]$ given by $\eta_M^\e(t) \coloneqq M(t\leq t+\e\one)$. Furthermore, we have the identities $M[\e][\e']= M[\e+\e']$ and $\eta_{M[\e']}^{\e} = \eta_M^\e[\e']$ for any $\e, \e'\geq 0$. 

  \begin{definition}\label{def:interleaving}
    Two persistence modules $M,N\colon \R^n\to\vect_{\field}$ are {\em $\e$-interleaved} if there exist morphisms $f\colon M\to N[\e]$ and $g\colon N\to M[\e]$ such that the following diagram commutes:
    \begin{equation}\label{eq:interleaving}
      \begin{gathered}
      \xymatrix@R=30pt@C=50pt{
      M \ar^-{\eta_M^\e}[r] \ar_<<<{f}|\hole[dr] & M[\e] \ar^-{\eta_M^\e[\e]}[r] \ar_<<<{f[\e]}|\hole[dr] & M[2\e] \\ 
      N \ar^-{\eta_N^\e}[r] \ar^<<<{g}[ur] & N[\e] \ar^-{\eta_N^\e[\e]}[r] \ar^<<<{g[\e]}[ur] & N[2\e] 
      }
      \end{gathered}
    \end{equation}
    Such a pair of morphisms~$(f,g)$ is called an {\em $\e$-interleaving} of the modules~$M,N$. The {\em interleaving distance} $\disti(M,N)$ is defined by:
    \[ \disti(M,N)\coloneqq \inf \left\{\e\geq 0 \mid \text{$M, N$ are $\e$-interleaved}\right\}. \]
  \end{definition}

When $\e=0$, the diagram~\eqref{eq:interleaving} collapses to the standard commutative diagram for isomorphism $M\simeq N$. In this sense, the interleaving distance measures the defect of isomorphism between $M$ and~$N$. It satisfies the triangle inequality, because interleavings compose to interleavings. It may take infinite values, when no interleaving between the modules exists. Finally, while  being isomorphic is equivalent to being $0$-interleaved and therefore implies being at interleaving distance~$0$, the converse is not true: for instance, some modules (e.g. $\field_{\{t\}}$, for any singleton $\{t\}\subset\R^n$) are called {\em ephemeral modules} because they are at interleaving distance~$0$ from the zero module without being isomorphic to it. Thus, the interleaving distance is an extended pseudo-distance on isomorphism classes of persistence modules. It has been shown to be an extended distance on the quotient of the persistence module category by the full subcategory of ephemeral modules~\cite{berkouk2021ephemeral,chazal2016observable}.

The next result justifies the use of the interleaving distance in the context of sublevel-sets topology. It says indeed that $\disti$ is a stable metric, in fact the most discriminative stable metric, in that context.

  \begin{proposition}
    \label{prop:inlearving_universality}
    For any topological space~$X$, functions $f,g\colon X\to\R^n$, and homology degree~$*$,
    \[ \disti(H_*(f),\, H_*(g)) \leq \|f-g\|_\infty \coloneqq \sup_{x\in X} \left\|f(x)-g(x)\right\|_\infty. \]
    Moreover, assuming either that $n=1$ or that $\field$ is a prime field, $\disti$ is the maximum extended pseudo-distance on~$\vect_{\field}^{\R^n}$ that satisfies the above inequality for any space~$X$, functions $f,g$, and homology degree~$*$.
  \end{proposition}
The second part of the proposition is due to Lesnick~\cite{lesnick}.  The first part follows immediately from the  commutativity of the diagram below for any  $\e\geq \|f-g\|_\infty$ and any $t\in\R^n$, which, by functoriality of~$H_*$, induces a commutative diagram like~\eqref{eq:interleaving} at the singular homology level, hence an $\e$-interleaving between $H_*(f)$ and $H_*(g)$: 
\[       \xymatrix@R=30pt@C=50pt{
      f^{-1}((-\infty, t]) \ar^-{\subseteq}[r] \ar_<<<{\subseteq}|\hole[dr] & f^{-1}((-\infty, t+\e\one]) \ar^-{\subseteq}[r] \ar_<<<{\subseteq}|\hole[dr] & f^{-1}((-\infty, t+2\e\one]) \\ 
      g^{-1}((-\infty, t]) \ar^-{\subseteq}[r] \ar^<<<{\subseteq}[ur] & g^{-1}((-\infty, t+\e\one]) \ar^-{\subseteq}[r] \ar^<<<{\subseteq}[ur] & g^{-1}((-\infty, t+2\e\one]) 
      } \]

    \subsection{Bottleneck distance}
  \label{sec:bottleneck}

  Definition~\ref{def:interleaving} does not provide any clear avenue for computing the interleaving distance. This is due to the size of the search space, which, up to shifts, is the product of the $\Hom$-spaces between the two modules under consideration. The idea behind the bottleneck distance is to turn the computation into a combinatorial assignment problem, whose search space is much smaller. For this, it leverages the direct-sum decompositions of the modules and restricts its focus to those interleavings that factor through these decompositions, considering partial matchings between the direct summands of the two modules under consideration.

  In the following, $M, N$ are two persistence modules over~$\R^n$, and $M\simeq \bigoplus_{j\in J} M_J$ and $N \simeq \bigoplus_{l\in L} N_l$ are their direct-sum decompositions as per Theorem~\ref{th:decomp}.

  \begin{definition}\label{def:bottleneck}
    A {\em partial matching} between the summands of $M$ and~$N$ is a bijection $J'\to L'$ between some subsets $J'\subseteq J$ and $L'\subseteq L$. The modules~$M, N$ are {\em bottleneck $\e$-interleaved} if there is a partial matching $J\supseteq J'\stackrel{\gamma}{\longrightarrow} L'\subseteq L$ such that:
    \[ \begin{cases}
    \disti\left(M_j,\, N_{\gamma(j)}\right) \leq \e & \forall j\in J',\\
    \disti\left(M_j,\, 0\right) \leq \e & \forall j \in J\setminus J',\\
    \disti\left(N_l,\, 0\right) \leq \e & \forall l \in L\setminus L'.
    \end{cases} \]
    We call such a partial matching a {\em bottleneck
      $\e$-interleaving} of the modules~$M,N$.  The {\em bottleneck
      distance} $\distb(M,\, N)$ is defined by:
    \[ \distb(M,\, N) \coloneqq \inf\left\{ \e \leq 0 \mid M, N\ \text{are bottleneck $\e$-interleaved}\right\}. \] 
  \end{definition}

  From a bottleneck $\e$-interleaving of~$M,N$, it is straightforward to build an $\e$-interleaving of~$M,N$, by taking the direct sum of all the interleaving morphisms between summands of $M$ and of $N$. By construction, this $\e$-interleaving factors through the direct-sum decompositions of $M$ and $N$ in that it restricts to $\e$-interleavings between their matched summands on the one hand, between their unmatched summands and the zero module on the other hand. Hence the name {\em bottleneck $\e$-interleaving}, and the following immediate inequality:
  \begin{equation}\label{eq:di_leq_db}
    \forall M,N\colon \R^n\to\vect_{\field},\quad \disti(M,\, N) \leq \distb(M,\, N).
  \end{equation}
  
  The converse inequality does not hold in general. Indeed, not every $\e$-interleaving of $M,N$ factors through their direct-sum decompositions, nor can be turned into one that does without increasing~$\e$. See e.g.~\cite{bjerkevik} for an example. However, the converse inequality does hold in the $1$-parameter setting, which yields the so-called {\em isometry theorem}.
  \begin{theorem}[Isometry \cite{bauer2015induced,chazal-silva-glisse-oudot,cohen2007stability,lesnick}]
    \label{th:isometry}
    For all $M,N\in \vect_{\field}^{\R}$, $\disti(M,\, N) = \distb(M,\, N)$.
  \end{theorem}

From a practical standpoint, this result, combined with Theorem~\ref{th:decomp_total-order}, implies that, for persistence modules over~$\R$, computing the bottleneck distance is a way of computing the interleaving distance~$\disti$ exactly and efficiently, since it boils down to a mere assignment problem between interval modules, whose  pairwise interleaving distances are easy to compute. For persistence modules over~$\R^n$ with $n>1$, by~\eqref{eq:di_leq_db} the bottleneck distance still provides an upper bound on~$\disti$, and in the special case of interval decomposable modules it is still computable in polynomial time~\cite{dey2019computing}. There also exist non-trivial lower bounds on~$\disti$ that can be computed in polynomial time, for instance the {\em matching distance}~\cite{kerber2018exact,landi}. However, computing the interleaving distance itself, whether exactly or within any factor less than~$3$, is known to be NP-hard over $\R^n$ when $n>1$, even for interval decomposable persistence modules~\cite{bjerkevik2020computing}.

\subsection{Signed bottleneck distance}
  \label{sec:signed_distances}

  Since the bottleneck distance~$\distb$ is computable and relevant on interval decomposable modules, and since signed barcodes are made of (pairs of multisets of) intervals, it is a natural idea to adapt $\distb$ to them.
Let~$\Int$ be a fixed set of intervals of~$\R^n$. Given any persistence modules~$M,N$ over~$\R^n$ that have well-defined signed barcodes $\barcode_\Int(M)$ and $\barcode_\Int(N)$ relative to~$\Int$, we define their {\em signed bottleneck distance} as follows:  
\begin{equation}\label{eq:signed_bottleneck}
  \hdistb^{\Int}(M,\, N) \coloneqq \distb\left(
  \field_{\barcode_\Int^+(M)} \oplus
  \field_{\barcode_\Int^-(N)},\ \field_{\barcode_\Int^+(N)} \oplus
  \field_{\barcode_\Int^-(M)} \right),
\end{equation}
where, for any multiset~$\barcode$ of intervals of~$\R^n$, $\field_{\barcode}$ denotes the interval decomposable module $\bigoplus_{\int\in \barcode} \field_\int$. The intuition behind this definition comes from geometric measure theory: if we view multisets of intervals as point measures on the space of intervals of~$\R^n$, then $\distb$ is the Wasserstein-$\infty$ distance for a specific choice of cost function on transport plans, and \eqref{eq:signed_bottleneck} is its signed version for signed measures. Note that the signed version does not satisfy the triangle inequality, however it satisfies a certain universality property that makes it canonical among the dissimilarities on signed barcodes~\cite{botnan2024bottleneck,oudot-scoccola}. 

Assume now that $\Int$ contains all the principal upsets of~$R^n$, and let $\ES_{\Int}$ be the class of short exact sequences on which the $\Hom$ functors $\Hom(\field_\int, -)$ are exact for all $\int\in\Int$.
%

  \begin{theorem}[Theorem~8.4 in~\cite{oudot-scoccola}]
    \label{th:stab_signed_barcodes}
Assume that:
\begin{itemize}
\item[(i)] $\left\{ t+\e\one \mid t\in\int \right\}\in\Int$ for all $\int\in\Int$ and $\e\geq 0$;
\item [(ii)] $m\coloneqq \mathrm{gldim}^{\ES_\Int}\!\!\left(\vect_{\field, \mathrm{fp}}^{\R^n}\right) <\infty$;
\item[(iii)] there is a constant $c\geq 0$ such that $\distb(\barcode(M),\, \barcode(N))\leq c\; \disti(M,\, N)$ for all $M, N\in\vect_{\field, \mathrm{fp}}^{\R^n}$ that are $\ES_\Int$-projective.
\end{itemize}
Then, for all $M,N\in\vect_{\field, \mathrm{fp}}^{\R^n}$:
\[ \hdistb^{\Int}\left(M,\, N)\right)
\leq c(m+1)\; \disti(M,N). \]
\end{theorem}

Note that assumption~(i) here differs from assumption~(a) in Theorem~8.4 of~\cite{oudot-scoccola}: it is in fact a sufficient condition for (a) to hold in our context, and one that is more readily verifiable. For instance, it clearly holds when  $\Int=\{\text{principal upsets}\}$ or $\Int=\hooks$.

For $\Int=\{\text{principal upsets}\}$, Hilbert's Syzygy Theorem provides~(ii) with $m=n$, while Bjerkevik's stability theorem for free modules~\cite{bjerkevik} provides~(iii) with $c=n-1$ (assuming $n>1$), hence we have $\hdistb^{\Int}\leq (n^2-1)\,\disti$ in this setting. For $\Int=\hooks$,  Theorem~\ref{th:gldim_Rn} provides~(ii) with $m=2n-2$, while an adaptation of Bjerkevik's stability proof to the case of hook modules (see~\cite{botnan2024bottleneck}) provides~(iii) with $c=2n-1$, hence we have $\hdistb^{\Int}\leq (2n-1)^2\,\disti$ in this setting.

  \begin{example}\label{ex:signed_bottleneck_distance}
    Let $M,N\colon: \R^2\to\vect_{\field}$ be the two interval modules whose supports and  minimal projective resolutions relative to hooks are depicted in Figures~\ref{fig:signed_bottleneck_distance_M} and~\ref{fig:signed_bottleneck_distance_N} respectively.   
\begin{figure}[tb]
  \centering
  \includegraphics[width=0.6\textwidth]{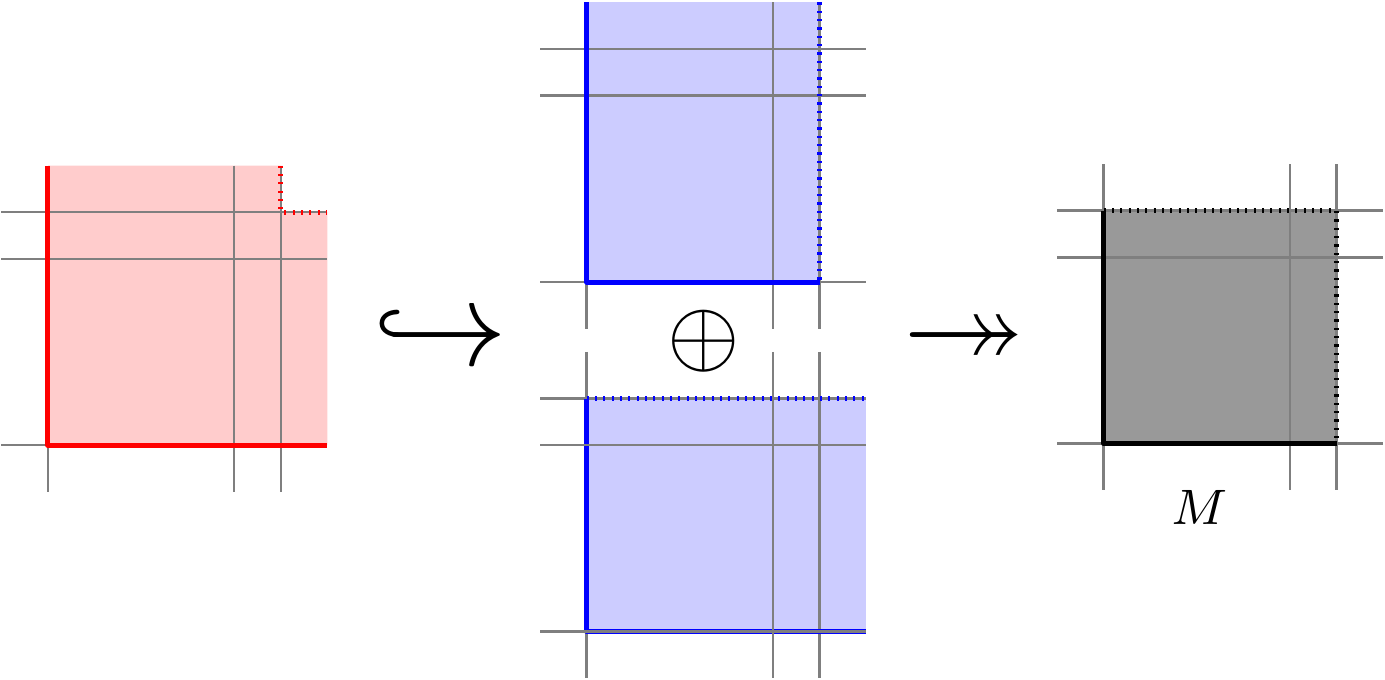}
  \caption{Minimal $\ES_{\hooks}$-projective resolution of the interval module~$M$ of Example~\ref{ex:signed_bottleneck_distance}.}
  \label{fig:signed_bottleneck_distance_M}
\end{figure}
\begin{figure}[tb]
  \centering
  \includegraphics[width=0.7\textwidth]{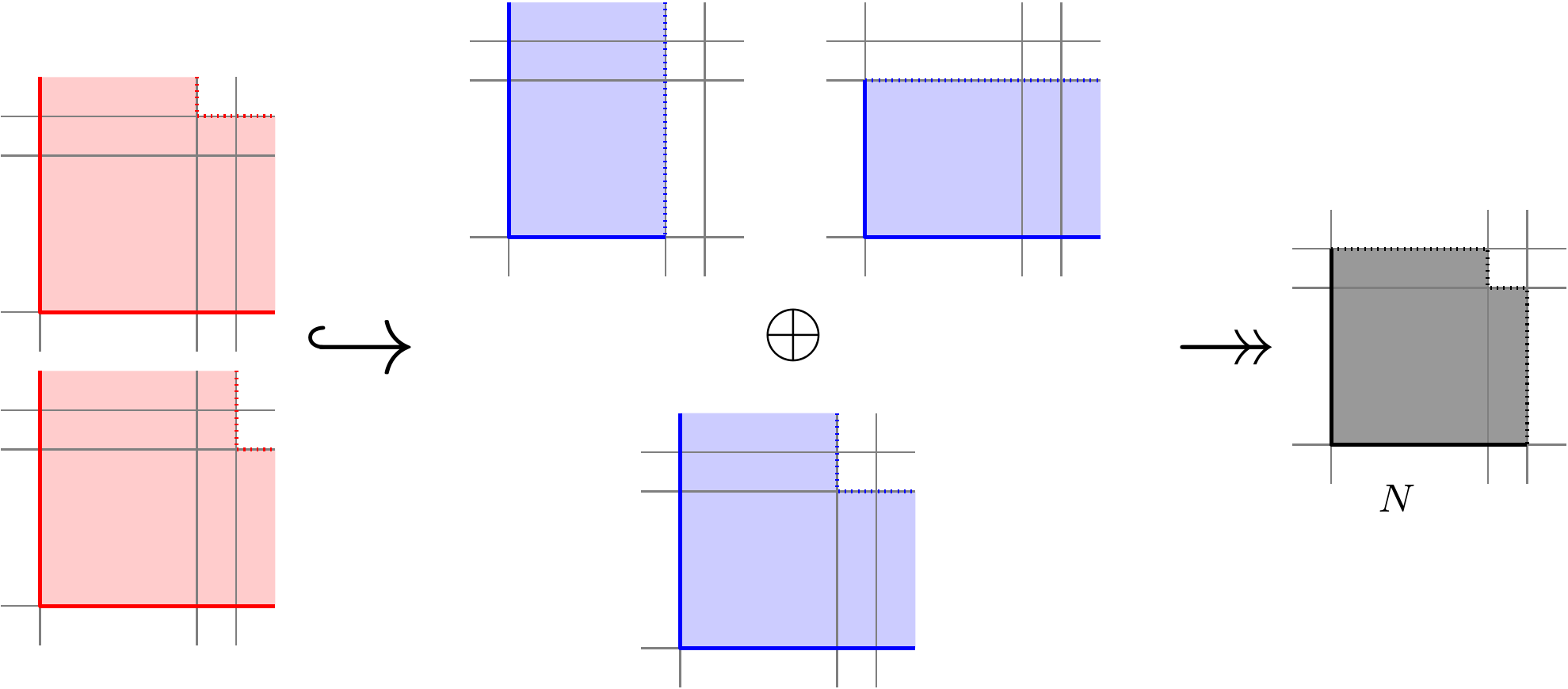}
  \caption{Minimal $\ES_{\hooks}$-projective resolution of the interval module~$N$ of Example~\ref{ex:signed_bottleneck_distance}.}
  \label{fig:signed_bottleneck_distance_N}
\end{figure}
Their signed bottleneck distance, realized by the partial matching between their signed barcodes shown in Figure~\ref{fig:signed_bottleneck_distance}, is equal to their interleaving distance, and so bounded from above by $(2\times 2-1)^2=9$ times their interleaving distance.   
\begin{figure}[tb]
  \centering
  \includegraphics[width=0.7\textwidth]{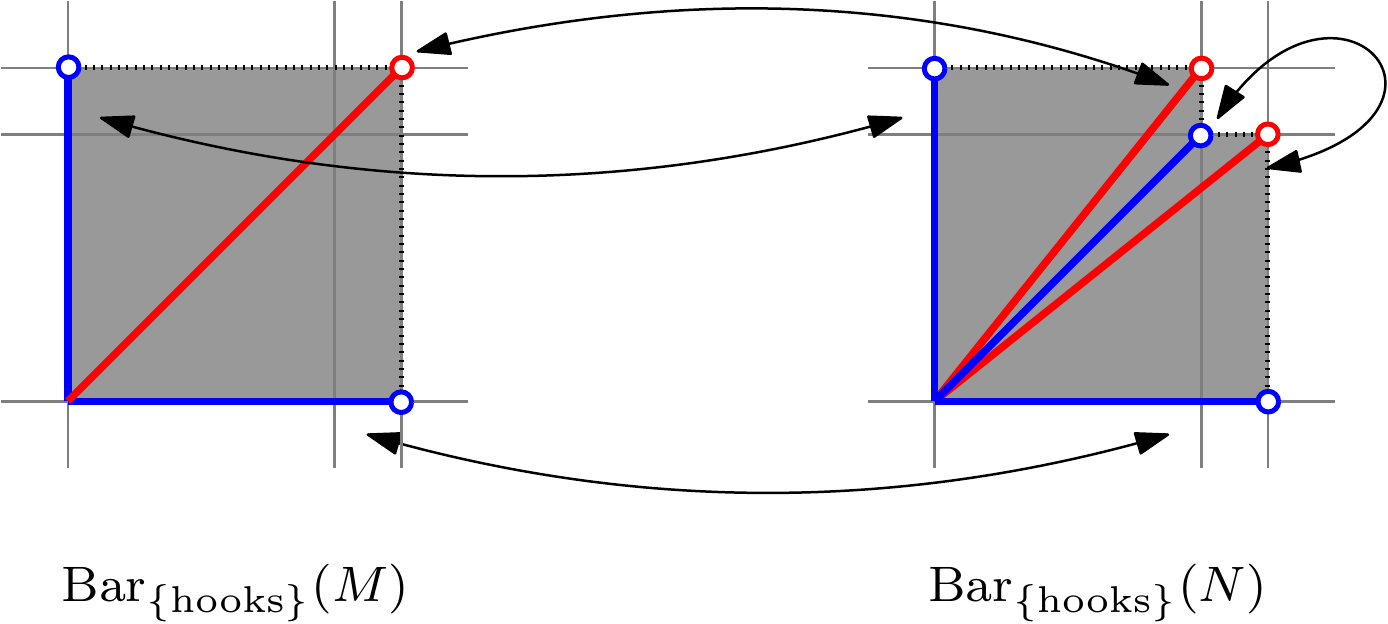}
  \caption{Partial matching realizing $\hdistb^{\hooks}(M,\, N)$ in Example~\ref{ex:signed_bottleneck_distance}.}
  \label{fig:signed_bottleneck_distance}
\end{figure}
  \end{example}

  \section{Differential calculus}
  \label{sec:diff_calc}

Let  $\R^n$ be our poset, equipped with the product order, and let $\vect_{\field, \mathrm{fp}}^{\R^n}$ be our persistence module category of interest. The finite presentedness assumption is motivated by the fact that, in TDA applications, persistence modules come from finite data. Note that $\vect_{\field, \mathrm{fp}}^{\R^n}$ is an essentially small category, so it makes sense to talk about set maps from or to it, provided we consider its objects up to isomorphism---which we do from now on. Our goal is to characterize what it means for a map from or to~$\vect_{\field, \mathrm{fp}}^{\R^n}$ to be differentiable, and to define derivatives for such maps so that the chain rule applies. The main obstruction is that the category does not come equipped with a differential structure on its object set, so we must define one.

\subsection{Setup}
\label{sec:diff_setup}
  
With TDA applications in mind, we place ourselves in the context of sublevel-sets topology, described in Section~\ref{sec:sublevel-sets}. Let $X$ be a fixed topological space, which we assume to be the underlying space of some finite simplicial complex~$K$. We consider functions $X\to\R^n$ that are defined simplexwise, meaning that they are the geometric realizations in~$X$ of functions $f=(f_1, \cdots, f_n)\colon K\to\R^n$. We further assume that every considered $f$ is {\em monotonic} in the sense that each $f_i$ is monotonic with respect to the face inclusion order on~$K$, that is: $f_i(\sigma)\leq f_i(\tau)$ for all simplices $\sigma\subseteq \tau\in K$.
These monotonic functions $K\to\R^n$ form a closed convex affine cone in $(\R^n)^{\# K}$, denoted~$(\R^n)^{\angle K}$,
which is the intersection of the closed half-spaces $\{f=(f_1,\cdots, f_n)\colon K\to \R^n \mid f_i(\sigma) \leq f_i(\tau)\}$ for all pairs of simplices $\sigma\subseteq \tau\in K$ and all coordinates $i\in \llbracket 1, n\rrbracket$. 
In the following, for the sake of simplicity of the exposition, we identify the functions $f\in (\R^n)^{\angle K}$ with their geometric realization in~$X$.

Every function $f\in (\R^n)^{\angle K}$ has a well-defined persistent (singular) homology $H_*(f)$ as a function $X\to\R^n$. Meanwhile, thanks to its monotonicity as a function $K\to\R^n$, $f$ also has a well-defined persistent (simplicial) homology, which is isomorphic to~$H_*(f)$. 
Since $K$ is finite, $H_*(f)$
is finitely presented and there is an upper bound, depending only on~$n$ and~$\# K$, on the size of its minimal presentation.  

In practice, the monotonic functions $K\to\R^n$ are built from data, so there is a map $\Phi\colon \R^m\to (\R^n)^{\angle K}$ from some parameter space that we assume to be the Euclidean space~$\R^m$ (or some subanalytic subset): for instance, $\R^m$ can be the space of labeled point clouds with $r$ points in $\R^d$ (in which case $m=dr$), or the parameter space of some neural network architecture. This setup is summarized in the first part of the following pipeline:

\[
\xymatrix@R=10pt@C=26pt{
\R^m  \ar_-{\Phi}[r] & \left(\R^n\right)^{\angle K} 
\ar_-{H_*(\cdot)}[r] 
\ar@/_2pc/_-{\Pers}[rr]|>>>>>>>\hole & \vect_{\field, \mathrm{fp}}^{\R^n}
\ar[r] 
\ar@/_2pc/_-{\loss}[rr]&
\barcodes  
  \ar_-{\Psi}[r] & \R
}
\]

The second part of the pipeline is for maps from the persistence module category. Generally speaking, one may consider maps to any Euclidean space (or any subanalytic subset), but anticipating the next section on optimization, we take $\R$ as the codomain. In order to leverage the theory developed in the previous sections, we assume that the map $\loss\colon \vect_{\field, \mathrm{fp}}^{\R^n} \to \R$ under consideration factors through some metric space of invariants, in our setting: either the space of finite barcodes equipped with the bottleneck distance (when $n=1$); or the space of finite signed barcodes equipped with the signed bottleneck distance (when $n>1$); in either case we write $\barcodes$ for the space of invariants. To complete the picture, we call {\em persistent homology map} and write~$\Pers$ for the composite map $(\R^n)^{\angle K}\to  \vect_{\field, \mathrm{fp}}^{\R^n} \to \barcodes$, and we write $\Psi$ for the map $\barcodes\to\R$ through which $\loss$ factors.

As the space of barcodes does not come equipped with a differential structure, we lift it to a disjoint union of Euclidean spaces as illustrated below:
\begin{equation} \label{eq:diff_pipeline_lift}
  \begin{gathered}
\xymatrix@R=10pt@C=26pt{
&&& \bigsqcup_{k=1}^\infty \R^{(2n+1)k} \ar@/^1ex/^-{q}[dd] \ar@{..>}[ddr] \\\\
\R^m  \ar_-{\Phi}[r] & \left(\R^n\right)^{\angle K} 
\ar@{..>}^-{~}[uurr]
\ar_-{H_*(\cdot)}[r] 
\ar@/_2pc/_-{\Pers}[rr]|>>>>>>>\hole & \vect_{\field, \mathrm{fp}}^{\R^n}
\ar[r] 
\ar@/_2pc/_-{\loss}[rr]&
\barcodes  
\ar@/^1ex/^-{p}[uu]  \ar_-{\Psi}[r] & \R
}
  \end{gathered}
\end{equation}
Specifically, we map each proper hook $[s,t)$ with $s<t\in\R^n$ to the vector $(s,t,1)\in \R^{2n+1}$, and each principal upset $[t, \infty)$ to the vector $(t,t,1)\in \R^{2n+1}$. Given a fixed total order on the bars (proper hooks and principal upsets) in $\R^n$, we now map each barcode with $k$ bars to the vector in $\R^{(2n+1)k}$ whose coordinates are given by the list of the coordinates of the images of its bars in $\R^{2n+1}$, sorted according to the total order on bars. In the following, we take the total order on bars induced by the lexicographic order on the coordinates of their images in $\R^{2n+1}$.  For example, in the case $n=1$, the barcode
    \[ \left\{[0, \infty),\, [-1, 3.2),\, [0, 2.5)\right\}\]
          is sent to the vector
          \[ (-1, 3.2, 1, 0, 0, 1, 0, 2.5, 1)\in\R^9; \]
          in the case $n=2$,  the barcode
          \[ \left\{[(0, 0), (3,2)),\, [(-1, 1), (2,2)),\, [(2,3), \infty)\right\} \]
                is sent to the vector
                \[ (-1, 1, 2, 2, 1, 0, 0, 3, 2, 1, 2, 3, 2, 3, 1)\in\R^{15}. \] 
                For signed barcodes $\barcode=(\barcode^+,\barcode^-)$, we map $\barcode^+$ to $\R^{(2n+1)\# \barcode^+}$ as above, and we map $\barcode^-$ to $\R^{(2n+1)\# \barcode^-}$ by sending each of its proper hooks $[s,t)$ with $s<t\in\R^n$ to $(s, t, -1)$ and each of its principal upsets $[t,\infty)$ to $(t, t, -1)$, keeping the same total order on bars as before. In other words, the last coordinate of the image of each bar encodes its sign in the signed barcode. Then, we concatenate the two vectors in $\R^{(2n+1)(\# \barcode^++ \# \barcode^-)}$. For instance, in the case $n=2$, the signed barcode
                    \[ \left( \left\{[(0,0) (0,1)),\, [(0,0), (1,1))\right\},\ \left\{ [(0,0), (1,0))\right\}\right)\]
is sent to the vector
                    \[ (0, 0, 0, 1, 1, 0, 0, 1, 1, 1, 0, 0, 1, 0, -1) \in \R^{15}.\]
                    Let us write $p\colon \barcodes\to\bigsqcup_{k=1}^\infty \R^{(2n+1)k}$ for the resulting lifting map on barcodes and signed barcodes. This map clearly admits a left inverse~$q$, since every (signed or regular) barcode can be reconstructed from its vectorization.
                    Note that $p$ is not onto, given the special form of the vectors obtained in $\R^{2n+1}$. Note also that $p$ is not continuous either, because barcodes with different numbers of bars are sent to different terms of the disjoint union $\bigsqcup_{k=1}^\infty \R^{(2n+1)k}$. However, $p$ is good enough for our purposes, as we will see in the next subsection.

                    As mentioned previously,
                    there is a global upper bound, depending only on $n$ and~$\# K$,  on the size of the minimal presentations of the persistent homology~$H_*(f)$ of the monotonic maps $f\colon K\to\R^n$. As a consequence, there is also a global upper bound on the size of their (signed or regular) barcode. Thus, in effect, $\barcodes$ is lifted  to a finite disjoint union of Euclidean spaces, which, in the following, we view as a subanalytic subset of some real analytic manifold.

\subsection{Results}

Proposition~\ref{prop:inlearving_universality}, combined with Theorem~\ref{th:isometry} (case $n=1$) or~\ref{th:stab_signed_barcodes} (case $n>1$), implies that the persistent homology map~$\Pers$ is Lipschitz continuous. In practice, the maps $\Phi$ and $\Psi$ are chosen by the user, and most common choices
are known to be locally Lipschitz continuous. Under this assumption, the composite map $\Psi\circ\Pers\circ\Phi\colon \R^m\to\R$ is also locally Lipschitz continuous, and  therefore differentiable almost everywhere in~$\R^m$, by Rademacher's theorem~\cite[Thm. 3.4.1 \& Rmk. 3.4.2]{cobzacs2019lipschitz}.
In the following we aim at getting some control over the singular set of the map in $\R^m$, and at deriving a chain rule to compute the map's derivatives at its regular points. For this we assume that the maps $\Phi$ and $\Psi\circ q$ are subanalytic\footnote{Strictly speaking, we see $\Psi\circ q$ as a map from the finite disjoint union of Euclidean spaces in which barcodes coming from monotonic maps $K\to\R^n$ are lifted, and as mentioned previously, we see this finite disjoint union as a subanalytic subset of some real analytic manifold.}, an assumption that is verified by most common choices of maps $\Phi$ and $\Psi$ in practice. The key point is that $p\circ \Pers$ has an elementary structure, described next. 
          
\begin{proposition}[\cite{carriere,scoccola-ICML-24,leygonie2022framework}]
  \label{prop:Pers_piecewise-affine}
There is an affine Whitney stratification of $(\R^n)^{\angle K}$ such that the map $p\circ\Pers$ is affine on each stratum.
\end{proposition}

In effect, the stratification is given by the restriction, to the cone~$(\R^n)^{\angle K}$, of the arrangement of hyperplanes $\{f=(f_1, \cdots, f_n)\colon K\to \R^n \mid f_i(\sigma) = f_i(\tau)\}$ for all pairs of simplices $\sigma \neq \tau\in K$ and all coordinates $i\in \llbracket 1, n \rrbracket$. The intuition behind the result is that, for any monotonic functions $f,g\colon K\to\R^n$ belonging to the same stratum, the coordinate functions $f_i$ and $g_i$ induce the same pre-order on the simplices of~$K$, so the persistent homologies $H_*(f)$ and $H_*(g)$ are isomorphic persistence modules up to a monotonic piecewise affine reparametrization of each coordinate axis in~$\R^n$; as a consequence, the (signed or regular) barcodes of $H_*(f)$ and $H_*(g)$ are the same up to the reparametrization of the coordinates of their bars' endpoints.

\begin{corollary}\label{cor:PsiPersPsi_definable}
Assume $\Phi$ and $\Psi\circ q$ are subanalytic. Then, $p\circ\Pers\circ\Phi$ and $\Psi\circ\Pers\circ\Phi$ are subanalytic as well, and there is a finite $C^1$ Whitney stratification of~$\R^m$ such that $p\circ\Pers\circ\Phi$ and $\Psi\circ\Pers\circ\Phi$ are differentiable on top-dimensional strata, $\Psi \circ q$ is differentiable on the images of these strata via $p\circ\Pers\circ\Phi$, and the following chain rule applies:
\[
\nabla_x \left(\Psi\circ \Pers\circ \Phi\right) = \nabla_{(p\circ\Pers\circ\Phi)(x)} \left(\Psi\circ q\right)\ 
\mathbf{J}_x \left(p \circ \Pers\circ \Phi\right).
\]
\end{corollary}

That $p\circ \Pers\circ \Phi$ is subanalytic follows from Proposition~\ref{prop:Pers_piecewise-affine} and the assumption that $\Phi$ itself is subanalytic. Then, the composition $\Psi\circ q\circ p \circ\Pers\circ\Phi=\Psi\circ\Pers\circ\Phi$ is also subanalytic. The existence of a finite $C^1$ Whitney stratification of~$\R^m$ such that both $p\circ \Pers\circ \Phi$ and $\Psi\circ\Pers\circ\Phi$ are differentiable on the top-dimensional strata and that $\Psi\circ q$ is differentiable on the images of these strats via $p\circ \Pers\circ \Phi$, follows from the Whitney stratifyability of subanalytic sets and maps~\cite{hardt1975topological,hardt1975stratification}. The chain rule in top-dimensional strata is then the usual chain rule for differentiable maps between differentiable manifolds.

\begin{remark}
Corollary~\ref{cor:PsiPersPsi_definable} extends to more general o-minimal structures for which Whitney stratifyability holds, such as for instance analytic-geometric categories~\cite{van1996geometric}. 
\end{remark}

\subsection{Examples}

We now explore two concrete examples of maps~$\Phi$ (Example~\ref{ex:Rips}) and~$\Psi$ (Example~\ref{ex:total_pers}) that are popular in applications. We provide explicit formulas for their derivatives, which can be combined together via the chain rule of Corollary~\ref{cor:PsiPersPsi_definable}.

\begin{example}[Rips filtration]
    \label{ex:Rips}
  Let $r\in \N^*$ be fixed, and let $K$ be the full simplex on~$\llbracket 1, r \rrbracket$, whose faces are in one-to-one correspondence with the power set $2^{\llbracket 1,r \rrbracket}\setminus\{\emptyset\}$. As before, $X$ is the underlying space of~$K$. For any finite set  $A=\{a_1, \cdots, a_r\}$ of labeled points in Euclidean space~$\R^d$, we consider the geometric realization in~$X$ of the following monotonic function $f_A\colon K\to\R$:
  \[
  \begin{cases}
    f_A(\{i_0\})=0 & \forall 1\leq i_0\leq r,\\
    f_A(\{i_0, \cdots, i_s\}) = \displaystyle \max_{0\leq j,l \leq s} \left\|a_{i_j} - a_{i_l}\right\|_2 & \forall s \geq 1\ \text{and}\ \forall 1\leq i_0< \cdots < i_s\leq r.
  \end{cases}
  \]
  This function is called the {\em Rips filtration} of the labeled point set~$A$. It serves as a combinatorial approximation to the distance to~$A$ in the ambient space~$\R^d$---see e.g. Chapters 4-5 in~\cite{oudot2015persistence} for the precise connection. Its parameters are the coordinates of the labeled points of~$A$ in~$\R^d$. We consider the subspace of monotonic functions $f_A\colon K\to\R$ for $A$ ranging over all the sets of $r$ labeled points in $\R^d$ (with possible repetitions of the points). Referring to the pipeline~\eqref{eq:diff_pipeline_lift}, we thus have $n=1$, $m=dr$, and $\Phi\colon A \mapsto f_A$. This map~$\Phi$ is known to be subanalytic~\cite{divol2019density}, and it admits a $C^1$ Whitney stratification where each top-dimensional stratum is the locus of the labeled point clouds that have nonzero and distinct pairwise distances inducing the same total order on the pairs of labels, while the union of the lower-dimensional strata is composed of those labeled point clouds with at least one trivial or two equal pairwise distance(s). The derivatives of~$\Phi$ at $A\in\R^{dr}$ in some top-dimensional stratum have the following elementary expression: for any simplex $\sigma\in K$,
  \[
  \left(\frac{\partial \Phi(A)}{\partial a_i}\right)_\sigma =
  \begin{cases}
    \frac{a_i - a_j}{\left\|a_i-a_j\right\|_2} & \text{if}\ \exists j\neq i \in \sigma \ \text{s.t.}\ \left\|a_i-a_j\right\|_2 = \displaystyle \max_{i',j'\in\sigma} \left\|a_{i'} - a_{j'}\right\|_2,\\
    0 & \text{otherwise}.
  \end{cases}
  \]
  Note that, even if $A$ itself is in a top-dimensional stratum~$S_A$, its image $f_A$ may be (in fact, is generally) located in a stratum~$T_A$ of the affine Whitney stratification of~$\R^{\angle K}$ that is not top-dimensional. This is because the appearance of an edge in the Rips filtration is concomitant to the appearance of all its cofaces that have it as their longest edge. Therefore, the map $p\circ\Pers$ per se may not be differentiable at~$f_A\in\R^{\angle K}$.

  The composite map $p\circ\Pers\circ\Phi$, however, is differentiable at~$A$. More precisely, the pairwise distances in each point cloud $A'\in S_A$ yield the same total order on the pairs of indices~$1\leq i<j\leq r$. In turn, this total order yields a unique pre-order~$\preceq_A$ on the simplices of~$K$ by the (increasing) values of~$f_A$---in other words, the point clouds in $S_A$ are mapped to the stratum~$T_A$. Then, fixing a priori a total order~$\leq_K$ on the simplices of~$K$, the point clouds in~$S_A$ induce the same total lexicographic order $\leq_A \coloneqq (\preceq_A,\, \leq_K)$.
  In turn, $\leq_A$ yields a unique assignment of the simplices of $K$ to the finite endpoints of the bars in the barcode of~$H_*(f_{A'})$ for any $A'\in S_A$, hence a unique assignment $\zeta$ of the simplices of $K$ to coordinates
  in the lifting space~$\R^{3k}$---where $k$ is the number of bars in the barcode of~$H_*(f_A)$. The derivatives of $p\circ \Pers \circ \Phi$ at~$A$ then have the following simple expression: for any coordinate $1\leq l \leq 3k$,
  \[
  \left(\frac{\partial \left(p\circ\Pers\circ\Phi\right)(A)}{\partial a_i}\right)_l =
  \begin{cases}
    \frac{a_i - a_j}{\left\|a_i-a_j\right\|_2} & \begin{minipage}{0.38\textwidth}if $l=\zeta(\sigma)$ for some $\sigma\in K$ and $\exists j\neq i \in \sigma$ such that\\$\left\|a_i-a_j\right\|_2 = \displaystyle \max_{i',j'\in\sigma} \left\|a_{i'} - a_{j'}\right\|_2,$\end{minipage}\\[2em]
    \left(\frac{\partial \left(p\circ\Pers\circ\Phi\right)(A)}{\partial a_i}\right)_{l-1} & \begin{minipage}{0.38\textwidth}if $l \equiv 2 \bmod 3$ and $(p\circ\Pers\circ\Phi)(A)_l = (p\circ\Pers\circ\Phi)(A)_{l-1}$ (bar endpoint at infinity),\end{minipage}\\[2em]
    0 & \text{otherwise}.
  \end{cases}
  \]
  Recall from our setup (Section~\ref{sec:diff_setup}) that $p$ sends infinite bars $[t,\infty)$ to vectors $(t,t,1)$, hence the second case in the above expression for the derivatives of $p\circ \Pers \circ \Phi$.
\end{example}

\begin{example}[Total persistence]
  \label{ex:total_pers}
  For $n=1$, let $\Psi\colon \barcodes\to\R$ map each finite barcode to the sum of the lengths of its finite bars, also known as its {\em total persistence}:
  \[ \forall \barcode\in\barcodes, \quad \Psi(\barcode) \coloneqq \sum_{\begin{smallmatrix}\int\in\barcode\\\sup I < \infty\end{smallmatrix}} (\sup I- \inf I). \] 
This sum is finite and can be rearranged as follows:
  \[ \forall \barcode\in\barcodes, \quad \Psi(\barcode) = \sum_{\begin{smallmatrix}\int\in\barcode\\\sup I < \infty\end{smallmatrix}} \sup I - \sum_{\begin{smallmatrix}\int\in\barcode\\\sup I < \infty\end{smallmatrix}} \inf I.  \]
  Pre-composing $\Psi$ by $q$ gives the following linear map $\R^{3k}\to\R$ for every barcode size~$k$:
\[ x \mapsto \sum_{i=0}^{k-1} x_{3i+2} - \sum_{i=0}^{k-1} x_{3i+1}. \]
One can indeed check that $\Psi\circ q$ as defined above yields $\Psi\circ q\circ p = \Psi$, the crux of the matter being that the contribution of every infinite bar~$I$ to $(\Psi\circ q\circ p)(\barcode)$ is zero as desired because $I$ is mapped to $(\inf I, \inf I, 1)$ by~$p$ then to $\inf I-\inf I = 0$ by $\Psi\circ q$.

The derivatives of $\Psi\circ q$ are then constant for every barcode size~$k$:
  \[ \frac{\partial (\Psi\circ q)}{\partial x_l} =
  \begin{cases}
    -1 & \text{if}\ l \equiv 1 \bmod 3,\\
    1 & \text{if}\ l \equiv 2 \bmod 3,\\
    0 & \text{otherwise}.
  \end{cases}
  \]
\end{example}

  \section{Optimization}
  \label{sec:optim}

  Referring again to the pipeline~\eqref{eq:diff_pipeline_lift}, we assume $\loss$ is some loss function on persistence modules that one seeks to minimize over the choice of parameter vector~$x\in\R^m$ from which the modules originate. We also assume the maps $\Phi$ and $\Psi\circ q$ have been chosen to be subanalytic and locally Lipschitz continuous---as most common choices in practice are. Then, by Corollary~\ref{cor:PsiPersPsi_definable} the composite map $F\coloneqq \Psi\circ\Pers\circ\Phi$ itself is subanalytic, so it is differentiable almost everywhere in~$\R^m$ and has a well-defined Clarke subdifferential $\partial F(x)$ at every point $x\in\R^m$:
  \[ \partial F(x)\coloneqq\Conv \left\{ \lim_{i \to \infty} \nabla F(y_i)  \mid  (y_i)_{i\in\N} \to x\ \text{and $F$ differentiable at every $y_i$}\right\}, \]
  where $\Conv$ is the notation for the convex hull in~$\R^m$. We consider the following differential inclusion:
  \[ \frac{dx}{dt} \in - \partial F(x(t))\ \text{for almost every}\ t\geq 0, \]
  whose solutions are the trajectories of the subgradient of~$F$.
  In order to approximate these trajectories, we use the stochastic subgradient descent algorithm, whose descent step writes as follows:
  \[ x_{i+1}: = x_i - \alpha_i ( g_i + \xi_i), \]
  where $(x_i)_{i\in\N}$ is the sequence of iterates, $(g_i\in \partial F(x_i))_{i\in\N}$ is any sequence of subgradients of~$F$ at the iterates, $(\alpha_i)_{i\in\N}$ is a sequence of non-negative learning rates, and $(\xi_i)_{i\in\N}$ is a sequence of random variables (the stochastic part of the descent). We make the following standard assumptions on these quantities---see Assumption~$C$ in~\cite{davis2020stochastic}:
  \begin{itemize}
  \item[(i)] $\sum_{i\in\N} \alpha_i = +\infty$ while $\sum_{i\in\N} \alpha_i^2 < +\infty$,
  \item[(ii)] $\sup_{i\in\N} \|x_i\|_2 < +\infty$ almost surely,
  \item[(iii)] almost surely, conditioned on the past, each $\xi_i$ has zero mean and its second moment is bounded from above by $\mu(x_i)$, for some fixed function  $\mu\colon \R^m\to\R_+$ that is bounded on bounded sets. 
  \end{itemize}

Under these assumptions, one can show that stochastic subgradient descent converges almost surely to critical points of~$F$.
    
  \begin{theorem}[Corollary~5.9 in \cite{davis2020stochastic}]
    \label{thm:SGD_conv}
    Given that $F$ is subanalytic and locally Lipschitz continuous, under assumptions (i)-(iii), almost surely every limit point~$x$ of the sequence of iterates $(x_i)_{i\in \N}$ is critical for~$F$ (i.e. $0\in\partial F(x)$) and the sequence of function values $(F(x_i))_{i\in\N}$ converges.
  \end{theorem}
  
We now give two illustrative examples of applications. For the sake of simplicity, we show basic applications that rely on the choices of maps $\Phi$ and $\Psi$ from Examples~\ref{ex:Rips} and~\ref{ex:total_pers}. Other, more involved applications have been considered in the literature, notably in deep learning contexts,  relying on different choices of maps $\Phi$ and $\Psi$---see e.g.~\cite{carriere,mukherjee2024d,scoccola-ICML-24} for some examples. 

\begin{example}[topological optimization]
  \label{ex:point_cloud_optim}
  Let $A$ be a finite set of points sampled iid from the uniform distribution on the square $[-1,1]^2\subset\R^2$, as illustrated in Figure~\ref{fig:point_cloud_optim}~(left).
  \begin{figure}[tb]
    \centering
    \includegraphics[width=0.3\textwidth]{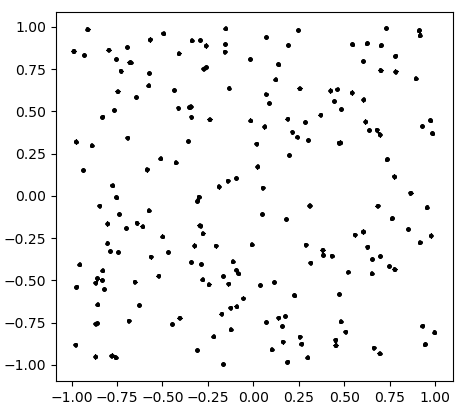}
    \hspace{0.03\textwidth}
    \includegraphics[width=0.3\textwidth]{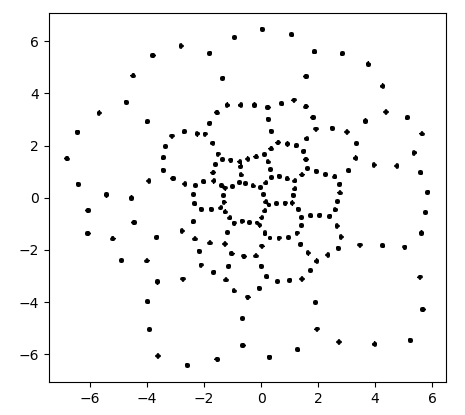}
    \hspace{0.03\textwidth}
    \includegraphics[width=0.3\textwidth]{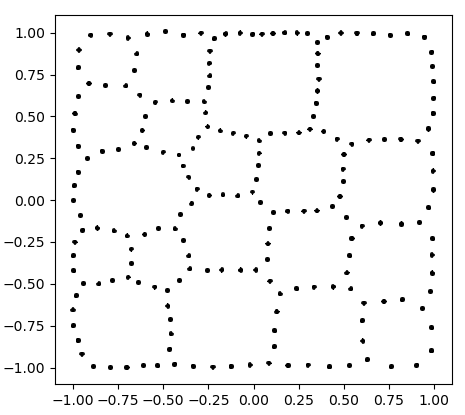}
    \caption{Left: the input point cloud~$A$. Center: the output after $100$ steps of stochastic subgradient descent on the functional $\Psi\circ\Pers\circ\Phi$, starting from~$A$, where persistent homology is taken in degree~$1$. Right: output after $100$ steps of stochastic subgradient descent on the regularized functional $(\Psi\circ\Pers\circ\Phi)(A) + \lambda\,\sum_{a\in A}\max\{0,\, \left\|a\right\|_\infty-1\}$.}
    \label{fig:point_cloud_optim}
  \end{figure}
  The goal is to move the points around in order to create as large `holes' as possible in the distribution. For this we use the setup of Examples~\ref{ex:Rips} and~\ref{ex:total_pers}, letting $d=2$, $n=1$, $m=2\,\# A$, $\Phi$ be the Rips filtration for labeled point clouds in~$\R^2$, and $\Psi$ be the opposite of the total persistence. Since we are interested in creating 1-dimensional holes, we consider homology in degree~$1$, so in the following $\Pers$ stands for persistent homology in degree~$1$. Our approach is to minimize the functional $\Psi\circ\Pers\circ\Phi$, that is, to maximize the total persistence in homology degree~$1$ of the Rips filtration of~$A$---recall that the Rips filtration of~$A$ serves as a combinatorial approximation to the distance to~$A$ in~$\R^2$, as mentioned in Example~\ref{ex:Rips}.

  Applying stochastic subgradient descent with standard parameter values to the functional $\Psi\circ\Pers\circ\Phi$, starting from~$A$,  gives the result shown in Figure~\ref{fig:point_cloud_optim}~(center) after $100$ steps. As expected, $1$-dimensional holes have started appearing in the distribution. However, the functional does not converge, as the optimization process keeps dispersing the points indefinitely. This is because scaling the data points by a factor $\nu>0$ results in scaling their pairwise distances by~$\nu$, and therefore scaling the functional by $\nu$ as well. And indeed, our initial objective of creating as large holes as possible in the distribution without any further requirements is ill-defined.  

  We therefore modify our objective by adding the constraint that the points of~$A$ should remain in the square $[-1,1]^2$. We treat this as a soft constraint and add the following regularization term to our functional, where $\lambda$ is a user-defined mixing parameter:
  \[  \lambda\,\sum_{a\in A}\max\{0,\, \left\|a\right\|_\infty-1\}. \] 
  Stochastic subgradient descent applied to the regularized functional eventually converges  thanks to the regularization term, and the output, depicted in Figure~\ref{fig:point_cloud_optim}~(right), qualitatively meets our new objective.
\end{example}

\begin{example}[topological regularization]
  \label{ex:autoencoder}
  Let $A$ be a finite set of data points in Euclidean space~$\R^3$, sampled from two entangled thickened circles, 
  as illustrated in Figure~\ref{fig:pt_autoencoder}~(left). The goal is to project the data down to $2$ dimensions, while preserving the geometry and underlying topological structure of the data as much as possible.
  \begin{figure}[tb]
    \centering
    \includegraphics[width=0.3\textwidth]{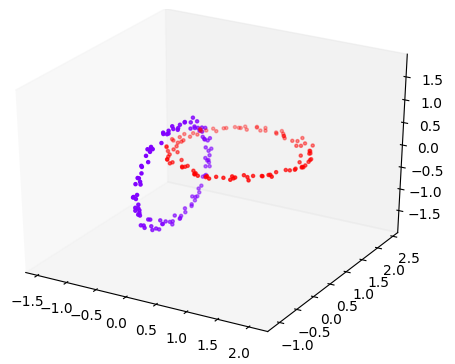}
    \hspace{0.03\textwidth}
    \includegraphics[width=0.3\textwidth]{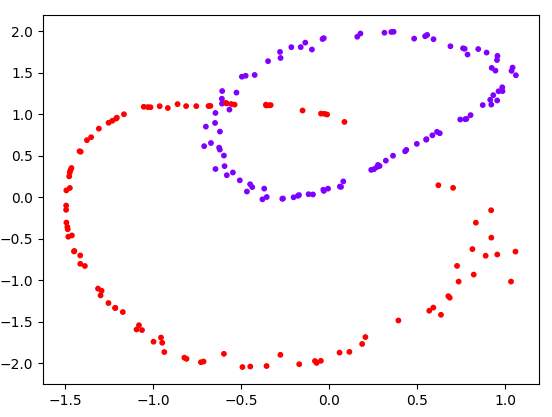}
    \hspace{0.03\textwidth}
    \includegraphics[width=0.3\textwidth]{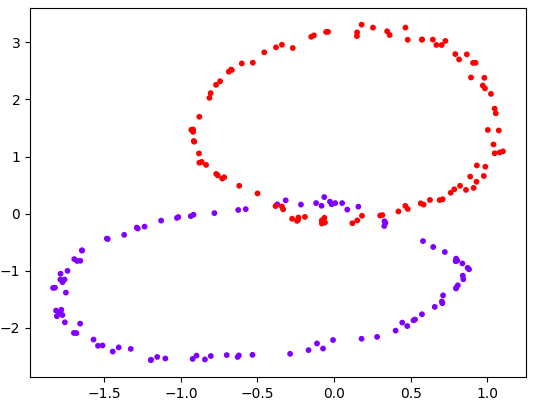}
    \caption{(from~\cite{carriere}) Left: the input point cloud~$A$ in $3$ dimensions, with data coming from the two thickened circles 
      marked in different colors. Center:  projection of~$A$ to $2$ dimensions by the autoencoder trained for $100$ epochs with the standard functional~\eqref{eq:autoenc_loss}. Right: projection of~$A$ to $2$ dimensions by the autoencoder trained for $100$ epochs with the regularized functional.}
    \label{fig:pt_autoencoder}
  \end{figure}

In order to do this, we use an auto-encoder architecture, which implements the composition of two maps:
\begin{itemize}
\item a projection map $\pi_\theta\colon \R^3\to\R^2$, called the {\em encoder}, parametrized by some vector $\theta\in\Theta$;
\item an embedding map $\iota_{\omega}\colon \R^2\to\R^3$, called the {\em decoder}, parametrized by some vector $\omega\in\Omega$.
\end{itemize}
The idea is then to make the composite map  $\iota_\omega \circ \pi_\theta\colon \R^3\to\R^3$ as close to the identity of~$\R^3$ as possible. For this, we set the parameters $\theta, \omega$ in such a way as to minimize the quadratic error between the points of~$A$ and their image through $\iota_\omega\circ\pi_\theta$:
\begin{equation}\label{eq:autoenc_loss}
  \min_{\theta, \omega}\ \sum_{a\in A} \left\|a-(\iota_\omega \circ \pi_\theta)(a)\right\|_2^2.
\end{equation}
The resulting projection of the point cloud~$A$ to $2$ dimensions after this training phase is shown in Figure~\ref{fig:pt_autoencoder}~(center). Beside the unavoidable metric distortion of the data, we see a break in the topological structure of its underlying two circles. In order to repair it, we add the following regularization term to our functional, which measures the discrepancy between the persistent homology in degree~$1$ of the Rips filtration of~$A$ and that of its projection~$\pi_\theta(A)$:
\[ \lambda\ \distone\left((\Pers\circ\Phi_3)(A),\, (\Pers\circ\Phi_2)(\pi_\theta(A))\right), \]
where $\lambda>0$ is a user-defined mixing parameter, $\Phi_d$ is the Rips filtration for labeled point clouds in~$\R^d$, $\Pers$ is the persistent homology map for homology in degree~$1$, and $\distone$ is a variant of the bottleneck distance between persistence modules (Definition~\ref{def:bottleneck}) where the maximum of the interleaving distances between their interval summands is replaced by the sum of these interleaving distances. As can be seen from  Figure~\ref{fig:pt_autoencoder}~(right), adding this regularization term to the functional helps the trained projection map~$\pi_\theta$ better preserve the topological structure of the two circles underlying the data. 
\end{example}

\section{Perspectives}
\label{sec:perspectives}

Many of the results presented in this text are preliminary and call for further investigations in order to consolidate and extend the parts of the theory of persistence modules they contribute to.
Here we report on a few concrete directions for future developments related to the differentiation, optimization, and application of the TDA pipeline.

Theorem~\ref{thm:SGD_conv}
provides convergence guarantees for stochastic subgradient descent on  functionals that factor through the persistence module category~$\vect_{\field, \mathrm{fp}}^{\R^n}$, but the result says nothing about the convergence rate. Theoretical bounds on the convergence rate were obtained in~\cite{leygonie2023gradient} for the case $n=1$, using a stratified variant of gradient sampling. Whether similar bounds can be obtained when $n>1$ remains open.

Strictly speaking, the framework described in Sections~\ref{sec:diff_calc} and~\ref{sec:optim} does not equip the persistence module category~$\vect_{\field, \mathrm{fp}}^{\R^n}$ itself with a differential structure, but rather takes a detour through  a certain lift of the space of (signed) barcodes. An important follow-up question is whether the differential structure on the lifting space can be pushed down to the space of (signed) barcodes, and from there, to the persistence module category~$\vect_{\field, \mathrm{fp}}^{\R^n}$ itself. This was shown in~\cite{leygonie2022framework} to be possible in the case $n=1$, using the theory of diffeological spaces~\cite{iglesias2013diffeology} and leveraging the Lipschitz continuity of the left inverse~$q$ to the lifting map~$p$ as well as the Isometry Theorem~\ref{th:isometry}. Whether this is also possible when $n>1$ remains open, and the main bottleneck is that, in this setting, (signed) barcodes are no longer a complete algebraic nor metric invariant of persistence modules.  

The question of computing (signed) barcodes efficiently is central for applications. While it is possible to use standard computer algebra systems like GAP~\cite{GAP4} for this purpose, in terms of algorithmic complexity these libraries do not scale up well with the size of the simplicial complexes on which filtrations are built---these simplicial complexes may contain millions of simplices in practice. Dedicated software libraries have been developped over the years by the TDA community, which include highly efficient routines to compute barcodes in the $1$-parameter setting~\cite{bauer2021ripser,hylton2019tuning,maria2014gudhi}. In the multi-parameter setting, effective algorithms to compute minimal projective resolutions of $2$-parameter persistence modules were recently proposed~\cite{kerber-rolle,lesnick-wright-2}, opening an avenue for computing their signed barcodes relative to principal upsets. However, computing minimal projective resolutions relative to other classes of projectives, or in higher numbers of parameters, remains largely open despite some promising developments~\cite{cacholski-guidolin-ren-scolamiero-tombari}.

Barcodes and their signed version are but one of many existing invariants for persistence modules. Indeed, given that barcodes are not a complete invariant in the multi-parameter setting, a wealth of other invariants have been proposed, among which:  the generalized persistence diagrams~\cite{kim-memoli}, elder-rule staircodes~\cite{cai2021elder}, multipersistence module approximations~\cite{loiseaux2022fast}, multigraded Hilbert series~\cite{harrington-otter-schenck-tillmann}, generalized rank invariant landscapes (GRIL)~\cite{xin2023gril}, meta-diagrams~\cite{clause2023meta}, or birth-death functions~\cite{mccleary2022edit}, to name a few. While some of them like GRIL are provably stable and differentiable~\cite{mukherjee2024d,xin2023gril}, many still have unknown theoretical behavior in terms of stability and differentiability. In order to prevent the parallel development of multiple specialized frameworks for differential calculus and optimization with each one of these invariants, there is a need for the emergence of a unified theory that would encompass most of them. A first attempt in this direction was made recently~\cite{scoccola-ICML-24}, encouraging further investigation.













\bibliographystyle{abbrv}
\bibliography{./biblio}


\end{document}